\documentclass[11pt]{amsart}

\makeatletter
\usepackage{amssymb}
\usepackage{latexsym}
\usepackage{amsbsy}
\usepackage{amsfonts}

\def\marginpar#1{\ignorespaces}

\newtheorem{theorem}[equation]{Theorem}
\newtheorem{proposition}[equation]{Proposition}
\newtheorem{lemma}[equation]{Lemma}
\newtheorem{corollary}[equation]{Corollary}

\theoremstyle{definition}
\newtheorem{remark}[equation]{Remark}

\newtheorem{example}[equation]{Example}

\numberwithin{equation}{section}

\def\AArm{\fam0 \rm}%
\newdimen\AAdi%
\newbox\AAbo%
\def\AAk#1#2{\setbox\AAbo=\hbox{#2}\AAdi=\wd\AAbo\kern#1\AAdi{}}%

\newcommand{\BBone}{{\ensuremath{{\AArm 1\AAk{-.8}{I}I}}}}

\def\eqref#1{(\ref{#1})}
\def\eqlabel#1{\def\@currentlabel{#1}}

\def\formula#1{\def\@tempa{#1}\let\@tempb\theequation\def\theequation{%
\hbox{#1}}\def\@currentlabel{(\theequation)}$$}
\def\endformula{\leqno\hbox{(\@tempa)}$$\@ignoretrue\let\theequation\@tempb}

\def\given{\hskip5\p@\relax\vrule\@width.4\p@\hskip5\p@\relax}

\newcommand{\open}[1]{%
\par\normalfont\topsep6\p@\@plus6\p@\trivlist\item[\hskip\labelsep\itshape#1%
\@addpunct{.}]\ignorespaces}

\DeclareRobustCommand{\close}[1]{%
  \ifmmode % if math mode, assume display: omit penalty etc.
  \else \leavevmode\unskip\penalty9999 \hbox{}\nobreak\hfill
  \fi
  \quad\hbox{$#1$}}

\newlength{\toskip}

\def \R {{\mathbb R}}
\def \Q {{\mathbb Q}}

\def \P {{\mathbb P}}

\def \N {{\mathbb N}}

\def \L {{\mathbb L}}

\def \phi {\varphi}

\def \Var {\textrm{Var}}
\def \Osc {\textrm{Osc}}
\def \Ent {\textrm{Ent}}
\def \Ipsi {{\mathcal I}_\psi}
\makeatother

\begin{document}
%%%%%%%%%%%%%%%%%%%%%%%%%%%%%%%%%%%%%%%%%%%%%%%%%%%%%%%%%%%%%%%%%%%%%%%
\date{\today}

\title[ TRENDS TO EQUILIBRIUM...]{TRENDS TO EQUILIBRIUM IN TOTAL VARIATION DISTANCE.}

%\author[F. Barthe]{\textbf{\quad {Franck} Barthe $^{\heartsuit}$ \, \, }}

%\address{{\bf {Franck} BARTHE},\\  Universit\'e Paul Sabatier
%Institut de Math\'ematiques. Laboratoire de Statistique et Probabilit\'es, UMR C 5583\\ 118 route
%de Narbonne, F-31062 Toulouse cedex 09.} \email{barthe@math.ups-tlse.fr}

\author[P. Cattiaux]{\textbf{\quad {Patrick} Cattiaux $^{\spadesuit}$ \, \, }}

\address{{\bf {Patrick} CATTIAUX},\\ Ecole Polytechnique, CMAP, F- 91128 Palaiseau cedex,
\\ and Universit\'e Paul Sabatier
Institut de Math\'ematiques. Laboratoire de Statistique et Probabilit\'es, UMR C 5583\\ 118 route
de Narbonne, F-31062 Toulouse cedex 09.} \email{cattiaux@cmapx.polytechnique.fr}

\author[A. Guillin]{\textbf{\quad {Arnaud} Guillin $^{\diamondsuit}$}}

\address{{\bf {Arnaud} GUILLIN},\\  Ecole Centrale Marseille and LATP
UMR CNRS 6632, \\Centre de Mathematiques et Informatique Technopôle Château-Gombert, \\39, rue F.
Joliot Curie, F-13453 Marseille Cedex 13.} \email{guillin@cmi.univ-mrs.fr}

%\author[C. Roberto]{\textbf{\quad {Cyril} Roberto $^{\clubsuit}$ \, \, }}

%\address{{\bf {Cyril}  ROBERTO},\\  Universit\'es
%de Marne-la-Vall\'ee et de Paris 12 Val-de-Marne. Laboratoire d'analyse et math\'ematiques
%appliqu\'ees, UMR 8050,\\
%   Boulevard Des\-cartes, Cit\'e Descartes, Champs
%sur Marne. F-77454 Marne-la-Vall\'ee cedex 2.} \email{roberto@math.univ-mlv.fr}

\maketitle
 \begin{center}

%\textsc{$^{\heartsuit}$ Universit\'e de Toulouse \quad and \quad Institut Universitaire de France}
%\medskip

 \textsc{$^{\spadesuit}$ Ecole Polytechnique \quad and \quad Universit\'e de Toulouse}
\medskip

\textsc{$^{\diamondsuit}$ Ecole Centrale de Marseille \quad and \quad Universit\'e de Provence}
\medskip

%\textsc{$^{\clubsuit}$ Universit\'e de Marne la Vall\'ee \quad and \quad Universit\'e Paris 12
%Val-de-Marne}
 \end{center}

\begin{abstract}
This paper presents different approaches, based on functional inequalities, to study the speed of
convergence in total variation distance of ergodic diffusion processes with initial law satisfying
a given integrability condition. To this end, we give a general upper bound ``\`a la Pinsker"
enabling us to study our problem firstly via usual functional inequalities (Poincar\'e inequality,
weak Poincar\'e,... ) and truncation procedure, and secondly through the introduction of new
functional inequalities $\Ipsi$. These $\Ipsi$-inequalities are characterized through
measure-capacity conditions and $F$-Sobolev inequalities. A direct study of the decay of Hellinger
distance is also proposed. Finally we show how a dynamic approach based on reversing the role of
the semi-group and the invariant measure can lead to interesting bounds.

\medskip

\noindent{\sc R\'esum\'e.} Nous \'etudions ici la vitesse de convergence, pour la distance en
variation totale, de diffusions ergodiques dont la loi initiale satisfait une int\'egrabilit\'e
donn\'ee. Nous pr\'esentons diff\'erentes approches bas\'ees sur l'utilisation d'in\'egalit\'es
fonctionnelles. La premi\`ere \'etape consiste \`a donner une borne g\'en\'erale \`a la Pinsker.
Cette borne permet alors d'utiliser, en les combinant \`a une procedure de troncature, des
in\'egalit\'es usuelles (telles Poincar\'e ou Poincar\'e faibles,...). Dans un deuxi\`eme temps
nous introduisons de nouvelles in\'egalit\'es appel\'ees $\Ipsi$ que nous caract\'erisons \`a
l'aide de condition de type capacit\'e-mesure et d'in\'egalit\'es de type $F$-Sobolev. Une \'etude
directe de la distance de Hellinger est \'egalement propos\'ee. Pour conclure, une approche
dynamique bas\'ee sur le renversement du r\^ole du semigroupe de diffusion et de la mesure
invariante permet d'obtenir de nouvelles bornes int\'eressantes.
\end{abstract}

\bigskip

\bigskip

\bigskip

\section{\bf Introduction, framework and first results.}\label{Intro}

We shall consider a dynamics given by a time continuous Markov process $(X_t, \P_x)$ admitting an
(unique) ergodic invariant measure $\mu$. We denote by $L$ the infinitesimal generator (and $D(L)$
the extended domain of the generator), by $P_t(x,.)$ the $\P_x$ law of $X_t$ and by $P_t$ (resp.
$P_t^*$) the associated semi-group (resp. the adjoint or dual semi-group), so that in particular
for any density of probability $h$ w.r.t. $\mu$, $\int P_t(x,.) h(x) \mu(dx) = P_t^* h \, d\mu$ is
the law of $X_t$ with initial distribution $h d\mu$. By abuse of notation we shall denote by
$P_t^* \nu$ the law of $X_t$ with initial distribution $\nu$.

Our goal is to describe the rate of convergence of $P_t^* \nu$ to $\mu$ in total variation
distance. Indeed, the total variation distance is one of the natural distance between probability
measures. If $d\nu = h d\mu$, this convergence reduces to the $\L^1(\mu)$ convergence.
\medskip

Trends to equilibrium is one of the most studied problem in various areas of Mathematics and
Physics. For the problem we are interested in, two families of methods have been developed during
the last thirty years.

The first one is based on Markov chains recurrence conditions (like the Doeblin condition) and
consists in finding some Lyapunov function. We refer to the works by Meyn and Tweedie
\cite{MT2,MT3,DMT} and the more recent \cite{Forro,Ver,DFG}. In a very recent work with D. Bakry
(\cite{BCG}), we have studied the relationship between this approach and the second one.

The second family of methods is using functional inequalities. It is this approach that we shall
follow here, pushing forward the method up to cover the largest possible framework. This approach relies mainly on the differentiation (with respect to time)  of a quantity like variance or entropy along the semigroup and a functional inequality enables then to use Gronwall's inequality to get the decay of the differentiated quantity. However, Due to the non
differentiability of the total variation distance, this direct method is no more possible. Let us then first give general upper bound on total
variation which will lead us to the relevant functional inequalities for our study.
\medskip

\subsection{\bf A general method for studying the total variation distance.}
\medskip

The starting point is the following elementary extension of the so called Pinsker inequality.
\begin{lemma}\label{lempinsker}
Let $\psi$ be a $C^2$ convex function defined on $\R^+$. Assume that $\psi$ is uniformly convex on
$[0,A]$ for each $A>0$, that $\psi(1)=0$ and that $\lim_{u \to +\infty} (\psi(u)/u) = +\infty$.
Then there exists some $c_\psi
> 0$ such that for all pair $(\P,\Q)$ of probability measures, $$\parallel \P - \Q
\parallel_{TV} \leq c_\psi \, \sqrt{I_\psi(\Q|\P)}\quad \textrm{ where } \quad I_\psi(\Q|\P)=\int
\psi\left(\frac{d\Q}{d\P}\right) d\P$$ if $\Q$ is absolutely continuous w.r.t. $\P$, and
$I_\psi(\Q|\P)=+\infty$ otherwise.
\end{lemma}
\begin{proof}
For $0\leq u \leq A$ it holds $$\psi(u) - \psi(1) - \psi'(1)(u-1) \geq \frac 12 \,
\left(\inf_{0\leq v \leq A} \psi''(v)\right) \, (u-1)^2 \, .$$ Thanks to convexity the left hand
side in the previous inequality is everywhere positive. Since $\lim_{u \to +\infty} (\psi(u)/u) =
+\infty$, it easily follows that there exists some constant $c$ such that for all $0\leq u$,
$$(u-1)^2 \, \leq \, c \, (1+u) \, \left(\psi(u) - \psi(1) - \psi'(1)(u-1)\right) \, .$$ Take the
square root of this inequality, apply it with $u=h(x)=(d\Q/d\P )(x)$, integrate w.r.t. $\P$ and
use Cauchy-Schwarz inequality. It yields $$\left(\int |h-1| d\P\right)^2 \leq c \, \left(\int
(1+h) d\P\right) \, \left(\int \big(\psi(h) - \psi(1) - \psi'(1)(h-1)\big) d\P\right) \, .$$ Since
$h$ is a density of probability the result follows with $c_\psi = \sqrt{2c}$.
\end{proof}
\begin{remark}\label{remlinear}
Note that we may replace the assumption $\psi(u)/u \to \infty$ by $\liminf_{u \to +\infty} \,
(\psi(u)/u) - \psi'(1) = d > 0$. For instance we may choose $\psi(u) = u - \frac 32 +
\frac{1}{u+1}$. \hfill $\diamondsuit$
\end{remark}
\medskip

The main idea now is to study the behavior of
\begin{equation}\label{eqIpsi}
t \, \mapsto \, I_\psi(t,h) = I_\psi(P^*_thd\mu|d\mu) = \int \psi(P_t^* h) \, d\mu
\end{equation}
as $t \to \infty$. Notice that with our assumptions $I_\psi(h)=I_\psi(0,h) \geq 0$ thanks to Jensen
inequality. To this end we shall make the following additional assumptions. The main additional
hypothesis we shall make is the existence of a ``carr\'e du champ'', that is we assume that there
is an algebra of uniformly continuous and bounded functions (containing constant functions) which
is a core for the generator and such that for $f$ and $g$ in this algebra
\begin{equation}\label{eqcarre}
L(fg) = f Lg + g Lf + \Gamma(f,g).
\end{equation}
We also replace $\Gamma(f,f)$ by $\Gamma(f)$.  Notice that with our choice there is a factor 2
which differs from many references, indeed if our generator is $\frac 12 \,  \Delta$,
$\Gamma(f)=|\nabla f|^2$ which corresponds to $L=\Delta$ in many references. The correspondence is
of course immediate changing our $t$ into $2t$.

We shall also assume that $\Gamma$ comes from a derivation, i.e. for $f$, $h$ and $g$ as before
\begin{equation}\label{eqderive}
\Gamma(fg,h) \, = \, f \Gamma(g,h) + g \Gamma(f,h) \, .
\end{equation}
The meaning of these assumptions in terms of the underlying stochastic process is explained in the
introduction of \cite{cat4}, to which the reader is referred for more details (also see
\cite{bakry} for the corresponding analytic considerations). Note that we may replace $L$ by $L^*$
without changing $\Gamma$.

Applying Ito's formula, we then get that for all smooth $\Psi$, and $f$ as before,
\begin{eqnarray}\label{eqchainrule}
L \Psi(f) & = & \frac{\partial \Psi}{\partial x}(f) \, Lf  \, + \, 1/2 \, \frac{\partial^2
\Psi}{\partial x^2}(f) \, \Gamma(f) \, ,\\ \Gamma(\Psi(f)) & = & (\Psi'(f))^2 \, \Gamma(f) \, .
\nonumber
\end{eqnarray}

Under these hypotheses, we immediately obtain
\begin{equation}\label{eqdecay}
\frac{d}{dt} \, I_\psi(t,h) \, = \, - \, \int \, 1/2 \, \psi''(P^*_t h) \, \Gamma(P^*_t h) \, d\mu
\, .
\end{equation}
It follows

\begin{proposition}\label{proptrendexpo}
There is an equivalence between
\begin{itemize}
\item for all density of probability $h$ such that $\int \psi(h) d\mu < +\infty$ , $$I_\psi(t,h)
\leq e^{-t/2C_\psi} \, I_\psi(h) \, ,$$ \item for all nice density of probability $h$,
\begin{equation}\label{eqineqpsi} \int \, \psi(h) \, d\mu \, \leq \, C_\psi \, \int \, \psi''(h)
\, \Gamma(h) \, d\mu \, .
\end{equation}
\end{itemize}
In this case the total variation distance $$\parallel P_t^*(h \mu) - \mu\parallel_{TV} \leq M_\psi
\, e^{-t/4C_\psi} \, I_\psi(h) \, ,$$ goes to 0 with an exponential rate.
\end{proposition}
When there exists $C_\psi$ such that for all nice $h$, (\ref{eqineqpsi}) holds for $\mu$, we will
say that $\mu$ verifies an ${\mathcal I}_\psi$-inequality. Note that the proof of this last
proposition is standard, the direct part is obtained by looking at $I_\psi(t,h) - I_\psi(h)$ when
$t$ goes to 0, while the converse part is a direct consequence of Gronwall lemma.

Slower decay can be obtained by weakening \eqref{eqineqpsi}. Indeed replace \eqref{eqineqpsi} by
\begin{equation}\label{eqineqpsiweak}
\int \, \psi(h) \, d\mu \, \leq \, \beta_\psi(s) \, \int \, \psi''(h) \, \Gamma(h) \, d\mu \, + s
\, G(h) \, ,
\end{equation}
supposed to be satisfied for all $s>0$ for some non-increasing $\beta_\psi$, and some real valued
functional $G$ such that $G(P_t^*h) \leq G(h)$. An application of Gronwall's lemma implies that
$$I_\psi(t,h) \leq \xi(t) \, \left(I_\psi(h)+G(h)\right)$$ with $\xi(t)=\inf \{s>0 \, , \,
2\beta_\psi(s) \, \log(1/s) \leq t \}$. Following R\"ockner-Wang \cite{RW01}, such an inequality
may  be called a weak $\Ipsi$-inequality. They consider the variance case, namely
$\psi(u)=(u-1)^2$, when the entropy case, namely $\psi(u)=u\log u$ is treated in \cite{CGG}. The
only known converse statement is in the variance case.
\smallskip

In this work we shall push forward this approach in order to give some rate of convergence for all
$h \in \L^1(\mu)$. The key is the following trick (see \cite{CGG} section 5.2): if $h\in
\L^1(\mu)$ and for $K>0$
\begin{eqnarray}\label{eqTV}
\int  |P^*_t h - 1| d\mu & \leq & \int  |P_t^* (h\wedge K) - P_t^* h| d\mu + \int  | P_t^*(h\wedge
K) - \int (h\wedge K) d\mu| d\mu + |\int (h\wedge K) d\mu - 1| \nonumber \\ & \leq & \int  |P_t^*
(h\wedge K) - \int (h\wedge K) d\mu| d\mu + 2 \int (h-K) \BBone_{h\geq K} d\mu
\end{eqnarray}
where we have used the fact that $P_t^*$ is a contraction in $L^1$. The second term in the right
hand sum is going to 0 when $K$ goes to $+\infty$, while the first term can be controlled by
$\sqrt{I_\psi(t,h\wedge K)}$ according to Lemma \ref{lempinsker}. More precisely, according to De
La Vall\'ee-Poussin theorem, $$\int h \, \phi(h) d\mu < +\infty$$ for some nonnegative function
$\phi$ growing to infinity. So
$$\int h \BBone_{h\geq K} d\mu \leq \frac{1}{\phi(K)} \, \left(\int h \, \phi(h) d\mu\right) \,
,$$ and we get, provided \eqref{eqineqpsiweak} is satisfied
\begin{equation}\label{eqtrend}
\int  |P^*_t h - 1| d\mu \, \leq \, c_\psi \, \sqrt{\xi(t) \, \left(I_\psi(h\wedge K)+G(h\wedge
K)\right)} + 2 \frac{1}{\phi(K)} \, \left(\int h \, \phi(h) d\mu\right) \, .
\end{equation}
\medskip

\subsection{\bf About this paper.}

Functional inequalities like \eqref{eqineqpsi} have a long story. When $\psi(u)$ behaves like
$u^2$ (resp. $u \, \log(u)$) at infinity, they are equivalent to the Poincar\'e inequality (resp.
the Gross logarithmic Sobolev inequality). We refer to \cite{logsob} for an introduction to this
topic. Many progresses in the understanding of such inequalities have been made recently. We refer
to \cite{BerZeg,RW01,BCR2,CGG} for their weak versions and to \cite{w00,w05,BCR1,BCR3,RZ06} for
the so called $F$-Sobolev inequalities. All these inequalities will be recalled and discussed
later. Links with long time behavior have been partly discussed in \cite{CatGui2,CGG,BCG}. Note
that in the recent \cite{RZ06}, the authors study the decay of $P_t f$ for $f$ belonging to
smaller spaces than $\L^2$.

Our aim here is to give the most complete description of the decay to 0 in total variation
distance using these inequalities, i.e. we want to give a general answer to the following question
: if a density of probability $h$ satisfies $\int\psi(h) d\mu < +\infty$ for some $\psi$ convex
at infinity, what can be expected for the decay to equilibrium in terms of a  functional
inequality satisfied by $\mu$ ?
\medskip

To better see what we mean, let us describe the contents of the paper.

In Section \ref{sectionineq} we recall old and recent results connected with Poincar\'e's like
inequalities and logarithmic Sobolev like inequalities. Recall that log-Sobolev is always stronger
than Poincar\'e. For short Poincar\'e (resp. log-Sobolev) inequality ensures an exponential decay
for densities such that $\int h^2 d\mu < +\infty$ (resp. $\int h \log h d\mu < + \infty$).
Actually we shall see in the examples of Section \ref{sectionineq} that these integrability
conditions can be replaced by $\int h^p d\mu < +\infty$ for some $p>1$ (resp. $\int h \log_+^\beta
h d\mu < + \infty$ for some $\beta
> 0$) with still an exponential decay. For less integrable densities, weak forms of Poincar\'e and
log-Sobolev inequalities furnish an explicit (but less than exponential) decay.
\smallskip

The questions are then :
\begin{itemize}
\item if $p>2$ and $\int |h|^p d\mu$ if finite, can we obtain some exponential decay with a weaker
functional inequality; \item if $u \log(u) \ll \psi(u) \ll u^2$, is it possible to characterize
$\Ipsi$-inequality, thus ensuring and exponential decay of $I(t,\psi)$; \item if $\psi(u) \ll u
\log(u)$ what can be said ? \end{itemize} \smallskip

The first question has a negative answer, at least in the reversible case, according to an
argument in \cite{RW01} (see Remark \ref{remarno}).
\smallskip

The answer to the second question is the aim of Section \ref{seccapacity}. It is shown that for
each such $\psi$ one can find a (minimal) $F$ such that exponential decay is ensured by the
corresponding $F$-Sobolev inequality (see \eqref{eqdeffsob} for the definition), and conversely
(see Theorem \ref{thmcapacity}, Theorem \ref{thmconvcapacity} and Remark \ref{remlemeilleur}).
These inequalities have been studied in \cite{w00,w05,BCR1,BCR3,RZ06}. A key tool here is the use
of capacity-measure inequalities introduced in \cite{BR03} and developed in
\cite{BCR1,BCR2,BCR3,CGG}. Hence for exponential decay we know how to interpolate between
Poincar\'e and log-Sobolev inequalities.
\smallskip

The third question is discussed in Section \ref{secTV}. This section contains essentially negative
results. A particular case is the ultracontractive situation, i.e. when $P_t h \in \L^2(\mu)$ for
all $h \in \L^1(\mu)$ and all $t>0$. Indeed if a weak Poincar\'e inequality is satisfied, the true
Poincar\'e inequality is also satisfied in this case, yielding a uniform exponential decay in
total variation distance. What we show in Section \ref{secTV} is that a direct study of the total
variation distance, or of the almost equivalent Hellinger distance, furnishes bad results i.e.
uniform (not necessarily exponential) decays are obtained under conditions implying
ultracontractivity.
\smallskip

The next Section \ref{secother} contains a discussion inspired by the final section of \cite{DLM},
namely, what happens if instead of looking at the density $P_t^* h$ with respect to $\mu$, one
looks at the density $1/P_t^* h$ with respect to $P_t^* h \, d\mu$, that is we look at
$d\mu/d\nu_t$ where $\nu_t$ is the law at time $t$. We show that a direct study leads to new
functional inequalities (one of them however is a weak version of the Moser-Trudinger inequality)
which imply a strong form of ultracontractivity (namely the capacity of all non-empty sets is
bounded from below by a positive constant). However, we also show that one can replace the
integrability condition on $h$ by a geometric condition ($h \mu$ satisfies some weak Poincar\'e
inequality) provided the Bakry-Emery condition is satisfied (see Theorem \ref{thmdelmo}). This
yields apparently better results than the one obtained in Section \ref{sectionineq} under the
log-Sobolev inequality (which is satisfied with the Bakry-Emery condition).

\section{\bf Examples using classical functional inequalities.}\label{sectionineq}

 In this section we shall show how to apply the general method in some classical cases.

\subsection{\bf Using Poincar\'e inequalities.}\label{secpoincare}

If we choose $\psi(u)=(u-1)^2$, \eqref{eqineqpsi} reduces to the renowned Poincar\'e inequality.
In this case Lemma \ref{lempinsker} reduces to Cauchy-Schwarz inequality. Recall what is obtained
in this case

\begin{theorem}\label{thmpoincare}
The following two statements are equivalent for some positive constant $C_P$
\begin{itemize}
\item[] \textbf{Exponential decay in $\L^2$.} \quad  For all $f \in \L^2(\mu)$, $$\parallel P_t f
- \int f d\mu
\parallel_2^2 \, \, \leq \, e^{- t/C_P} \,  \parallel f - \int f d\mu \parallel_2^2 \, .$$
\item[] \textbf{Poincar\'e inequality.} \quad For all $f \in D_2(L)$ (the domain of the Friedrichs
extension of $L$) , $$\Var_\mu(f) \, := \, \parallel f - \int f d\mu \parallel_2^2 \, \, \leq \,
C_P \, \int \, \Gamma(f) \, d\mu \, .$$
\end{itemize}
Hence if a Poincar\'e inequality holds, for $\nu = h \mu$ with $h\in \L^2(\mu)$, $$\parallel
P_t^*\nu - \mu
\parallel_{TV} =
\parallel P_t^*h - 1\parallel_{\L^1(\mu)} \leq e^{- t/2 C_P} \, \parallel h -
1\parallel_{\L^2(\mu)}\, .$$
\end{theorem}
\smallskip

\begin{corollary}\label{corpoincare}
Let $\tilde{\phi}(u)=\sqrt{u} \phi(u)$ and $\tilde{\phi}^{-1}$ its inverse. If a Poincar\'e
inequality holds,
\begin{equation}\label{eq1}
\parallel P_t^*h - 1\parallel_{\L^1(\mu)} \leq \frac{4 \, \int h \, \phi(h) d\mu}
{(\phi\circ \tilde{\phi}^{-1})(2 \, \left(\int h \, \phi(h) d\mu\right) \, e^{t/2 C_P})} \, .
\end{equation}
\end{corollary}
\begin{proof}
First remark that
$$\Var_\mu(h\wedge K) \leq \int (h\wedge K)^2 \, d\mu \leq K \, \int (h\wedge K) d\mu \leq K \,
.$$ We may now use Cauchy-Schwarz inequality to control $\int  |P_t^* (h\wedge K) - \int (h\wedge
K) d\mu| d\mu $ by the square root of the Variance in \eqref{eqTV}. The result then follows by an
easy optimization in $K$. More precisely we may choose $K$ in such a way that both terms in the
right hand side of \eqref{eqTV} are equal.
\end{proof}
\smallskip

\begin{example}\label{examplepower}
Assume that $h \in \L^q(\mu)$ for some $1<q<2$. If a Poincar\'e inequality holds, \eqref{eq1}
yields after some elementary calculation,
\begin{equation}\label{eq2}
\parallel P_t^*h - 1\parallel_{\L^1(\mu)} \leq 4^{\frac{q}{2q-1}} \left(\int h^q d\mu\right)^{\frac{1}{2q-1}} \,
e^{- \, \frac{(q-1) t}{(2q-1) \, C_P}} \, .
\end{equation}
Note that for $q=2$ we do not recover the good rate $e^{-t/2C_P}$ but $e^{-t/3C_P}$. It is however
not surprising, the truncation method is robust but not so precise.

Note that, up to the constants, a similar result already appeared in \cite{w04}. Indeed if a
Poincar\'e inequality holds then for $1 \leq p < 2$, $$\int f^2 d\mu - \left(\int f^p
d\mu\right)^{2/p} \, \leq \, C_P \, \int \, \Gamma(f) d\mu$$ which is a Beckner type inequality,
called ($I_p$) in \cite{w04}. According to Proposition 4.1 in \cite{w04} (recall the extra factor
2 therein), $$\int (P_t^*h)^{2/p} d\mu - 1 \, \leq \, e^{- \, \frac{(2-p) t}{C_P}} \, \left(\int
h^{2/p} d\mu - 1\right) \, ,$$ so that, taking $p=2/q$ and applying Lemma \ref{lempinsker} with
$\psi(u)=u^q - 1$ in the left hand side of the previous inequality we recover an exponential rate
of decay as in \eqref{eq2}, but this time with the good constant in the exponential term. It is
once again a motivation to study $\Ipsi$-inequality.\hfill $\diamondsuit$
\end{example}
\smallskip

Remark that in the derivation of the Corollary we only used Poincar\'e's inequality for bounded
functions. Hence we may replace it by its weak form introduced in \cite{BerZeg,RW01}, that is, we
take for $G$ the square of the Oscillation of $h$ in \eqref{eqineqpsiweak}. It yields

\begin{theorem}\label{thmweakpoincare}\textbf{(\cite{RW01} Theorem 2.1)} \quad
Assume that there exists some non-increasing function $\beta_{WP}$ defined on $(0,+\infty)$ such
that for all $s>0$ and all bounded $f\in D_2(L)$ the following inequality holds
$$\textbf{Weak Poincar\'e inequality.} \qquad \Var_\mu(f) \, \leq \,
\beta_{WP}(s) \, \int \, \Gamma(f) \, d\mu \, + \, s \, \Osc^2(f) \, .$$ Then $$\Var_\mu(P_t^* f)
\leq 2 \xi_{WP}(t) \, \Osc^2(f) \quad \textrm{ where } \quad \xi_{WP}(t) = \inf \, \{s>0 \, , \,
\beta_{WP}(s) \, \log(1/s) \, \leq \, t\} \, .$$

Hence if a weak Poincar\'e inequality holds, $$\parallel P_t^*h - 1\parallel_{\L^1(\mu)} \leq
\frac{ 4 \, \int h \, \phi(h) d\mu}{(\phi \circ \theta^{-1})(\sqrt 2 \, \left(\int h \, \phi(h)
d\mu\right) /\sqrt{\xi_{WP}(t)})}  \, ,$$ where $\theta(u)= u \phi(u)$.
\end{theorem}
The proof of the last statement is similar to the proof of \eqref{eq1}.
\smallskip

\begin{remark}\label{remzitt}
Since we are interested in functions such that $\int h \, \phi(h) \, d\mu < +\infty$, instead of
using the truncation argument we may directly try to obtain a weak inequality with $G(h)=\parallel
h - m_h\parallel_\zeta$ where $\parallel. \parallel_\zeta$ denotes the Orlicz norm associated to
$\zeta(u)=u \phi(u)$, and $m_h$ is a median of $h$. Actually as shown in \cite{Zitt} Theorem 29,
provided $\phi(h)\geq h$, such an inequality is equivalent to the weak Poincar\'e inequality
replacing $\beta_{WP}(s)$ by
\begin{equation}\label{eqzitt}
\beta_\zeta(s)= 6 \, \beta_{WP}\left(\frac 14 \, \bar{\zeta}(s/2)\right) \, \textrm{ where }
\bar{\zeta}(u)=\frac{1}{\gamma^*(1/u)} \textrm{ with } \gamma(u)=\zeta(\sqrt u) \, ,
\end{equation}
and $\gamma^*$ is the Legendre conjugate of $\gamma$ (assumed to be a Young function here). Since
for a density of probability $m_h\leq 2$ and since there exists a constant $c$ such that
$$\parallel g \parallel_\zeta \leq c \, \left(1+\int g \, \phi(g) \, d\mu\right) \, ,$$
at least if $\zeta$ is moderate, we immediately get a decay result
\begin{equation}\label{eqzitt2}
\parallel P_t^*h - 1\parallel_{\L^1(\mu)} \leq  C \, \sqrt{\xi_\zeta(t)} \, \int h \, \phi(h) \, d\mu \,
,
\end{equation}
with $$\xi_\zeta(t) \, = \, \inf \, \{s \, ; \, \beta_\zeta(s) \, \log(1/s) \, \leq \, t \} \, .$$
If $\phi(u)=u^{p-1}$ for some $p>1$, Theorem \ref{thmweakpoincare} yields a rate of decay
$(\xi_{WP}(t))^{\frac{p-1}{2p}}$.

Similarly, but if $p>2$, up to the constants, $\gamma(u)=u^{p/2}$, $\gamma^*(u)= u^{p/(p-2)}$
hence $\bar{\zeta}(u)=u^{p/(p-2)}$ so that we get $\xi_\zeta(t)=(\xi_{WP}(pt/(p-2)))^{(p-2)/p}$
hence a worse rate of decay. \hfill $\diamondsuit$
\end{remark}
\medskip

Of course our approach based on truncation extends to many other situations, in particular if we
assume that $\int h \log h d\mu < +\infty$, a Poincar\'e inequality yields a polynomial behavior
$$\parallel P_t^*h - 1\parallel_{\L^1(\mu)} \leq C\left(\int h \log h d\mu\right) \, (C_P/t) \,
.$$ It was shown in \cite{CatGui1} that a Poincar\'e inequality is equivalent to a restricted
logarithmic Sobolev inequality (restricted to bounded functions). The truncation approach together
with this restricted inequality do not furnish a better result. However with some extra
conditions, which are natural for diffusion processes on $\R^n$, one can prove sub-exponential
decay. We refer to \cite{CGG} sections 4 and 5 for a detailed discussion.
\medskip

\subsection{\bf Using a logarithmic Sobolev inequality.}\label{seclogsob}

In the previous subsection we have seen (Example \ref{examplepower}) that a Poincar\'e inequality
implies an exponential decay for the total variation distance, as soon as $\nu = h \mu$ for $h\in
\L^q(\mu)$ for some $q>1$. In this section we shall see that a similar result holds if $\int h
\log_+^\beta h d\mu < +\infty$ for some $\beta > 0$, as soon as a logarithmic Sobolev inequality
holds. First of all we recall the following (corresponding to $\psi(u)=u \log u$ in the
introduction)

\begin{theorem}\label{thmlogsob}
The following two statements are equivalent for some positive constant $C_{LS}$
\begin{itemize}
\item[] \textbf{Exponential decay for the entropy.} \quad  For all density of probability $h$
$$\int  P_t^* h \log (P_t^* h) d\mu \, \leq \, e^{- 2t/C_{LS}} \,  \int h \log h d\mu \, .$$
\item[] \textbf{Logarithmic Sobolev inequality.} \quad For all $f \in D_2(L)$  , $$\Ent_\mu(f^2)
\, := \, \int f^2 \, \log \left(\frac{f^2}{\int f^2 d\mu}\right) d\mu \, \leq \, C_{LS} \, \int \,
\Gamma(f) \, d\mu \, .$$
\end{itemize}
Hence if a logarithmic Sobolev inequality holds, for $\nu = h \mu$ with $\Ent_\mu(h)<+\infty$,
$$\parallel P_t^*\nu - \mu
\parallel_{TV} =
\parallel P_t^*h - 1\parallel_{\L^1(\mu)} \leq e^{- t/C_{LS}} \, \sqrt{2 \Ent_\mu(h)} \, .$$
\end{theorem}

\begin{corollary}\label{corlogsob}
Define $\bar{\phi}(u)= \phi(u) \, \sqrt{\log u}$ for $u\geq 1$. Then if a logarithmic Sobolev
inequality holds,
\begin{equation}\label{eq3}
\parallel P_t^*h - 1\parallel_{\L^1(\mu)} \leq \frac{4 \, \int h \, \phi(h) d\mu}{(\phi \circ
\bar{\phi}^{-1})\left( \left(\int h \, \phi(h) d\mu\right) \, e^{t/C_{LS}}\right)} \, .
\end{equation}
\end{corollary}
\begin{proof}
The proof is similar to the one of Corollary \ref{corpoincare}, replacing the Variance by the
Entropy, Cauchy-Schwarz inequality by Pinsker inequality an using the elementary
$$\Ent_\mu(h\wedge K) \leq \int (h\wedge K) \log (h\wedge K) d\mu + \frac 1e \leq \log K + \frac
1e \, .$$ We may then assume that $K>e^{1/e}$ and make an optimization in $K$.
\end{proof}

\smallskip

\begin{example}\label{examplepentropy}
Assume that $\int h \log_+^\beta h \, d\mu < +\infty$ for some $0<\beta \leq 1$. The previous
result yields, provided a logarithmic Sobolev inequality holds,
\begin{equation}\label{eq4}
\parallel P_t^*h - 1\parallel_{\L^1(\mu)} \leq \left(\int h \log_+^\beta h \,  d\mu\right)^{\frac{1}{2\beta+1}} \,
 e^{-2\beta t/(2\beta+1) C_{LS}} \, .
\end{equation}
Hence here again we get an exponential decay provided some ``$\beta$-entropy'' is finite.

Actually, as in Example \ref{examplepower}, a similar result can be obtained, provided a
log-Sobolev inequality holds, using the more adapted
$$I_\psi(t,h) \, \leq \, e^{- \, C \beta t} \, I_\psi(h)$$ with $\psi(u)= u \, (\log^\beta(2+u) -
\log^\beta(3)) = u F(u)$. It will be the purpose of the next section. In fact, in this example,
we will even show that the assumption of a logarithmic Sobolev inequality to hold is not
necessary, a well adapted $F$-Sobolev inequality will be sufficient.
%The proof of this result can be done by showing that, if a log-Sobolev
%inequality holds, $$\int h \, F(h) d\mu \leq C \, \int \psi''(h) \, \Gamma(h) \, d\mu \, ,$$ using
%the measure-capacity characterization of the log-Sobolev inequality. We refer to \cite{BCR1}
%section 5.4, Theorem 20 and Lemma 21. In particular the claimed result can be shown by mimicking
%the proof of Theorem 20 therein with $h=f^2$. Also see the next section.
 \hfill $\diamondsuit$
\end{example}

\begin{remark}
It is well known that a logarithmic Sobolev inequality implies a Poincar\'e inequality. Hence we
may ask whether some stronger inequality than the log-Sobolev inequality, furnishes some
exponential decay under weaker integrability conditions. But here we have to face a new problem :
indeed classical stronger inequalities usually imply that $P_t$ is ultracontractive (i.e. maps
continuously $\L^1(\mu)$ into $\L^\infty(\mu)$). Hence in this case we get an exponential decay
for the $\L^1(\mu)$ norm, combining ultracontractivity and Poincar\'e inequality for instance. We
shall give some new insights on this in one of the next sections.

Examples of ultracontractive semi-groups can be found in \cite{Dav,KKR}. \hfill $\diamondsuit$
\end{remark}
\medskip

\begin{remark}
Since a logarithmic Sobolev inequality is stronger than a Poincar\'e inequality, it is interesting
to interpolate between both inequalities. Several possible interpolations have been proposed in
the literature, starting with \cite{LO00}. In \cite{BCR1} a systematic study of this kind of
$F$-Sobolev  inequalities is done. Note that a homogeneous $F$-Sobolev  inequality is written as
$$\int f^2 \, F\left(\frac{f^2}{\int f^2 d\mu}\right) \, d\mu \, \leq \, \int \, \Gamma(f) \, d\mu \, $$
hence does not correspond to \eqref{eqineqpsi}. That is why such inequalities are well suited for
studying the convergence of $P_t f$ (see \cite{RZ06}), while we are interested here in the
convergence of $P_t (f^2)$. Moreover their convergence are stated in Orlicz norm  (clearly adapted
to $F$-Sobolev), whereas ours are in more usual integral form.

The case of $F=\log$ corresponding to the log-Sobolev (or Gross) inequality appears as a very
peculiar one since it is the only one for which the $F$-Sobolev inequality corresponds exactly to
\eqref{eqineqpsi}. It is thus natural to expect that the weak logarithmic Sobolev inequalities are
well suited to furnish a good interpolation scale between Poincar\'e and Gross inequalities. This
point of view is developed in \cite{CGG}. We shall recall and extend some of these results below.
\hfill $\diamondsuit$
\end{remark}

Here again we may replace the logarithmic Sobolev inequality by a weak logarithmic Sobolev
inequality

\begin{theorem}\label{thmweaklogsob}\textbf{(\cite{CGG} Proposition 4.1)} \quad
Assume that there exists some non-increasing function $\beta_{WLS}$ defined on $(0,+\infty)$ such
that for all $s>0$ and all bounded $f\in D_2(L)$ the following inequality holds
$$\textbf{Weak log-Sobolev inequality.} \qquad \Ent_\mu(f^2) \, \leq \,
\beta_{WLS}(s) \, \int \, \Gamma(f) \, d\mu \, + \, s \, \Osc^2(f) \, .$$ Then for all
$\varepsilon >0$, $\Ent_\mu(P_t^* h) \leq  (\frac 1e + \varepsilon) \, \xi_{WLS}(\varepsilon,t) \,
\Osc^2(\sqrt{h})$  where $\xi_{WLS}(\varepsilon,t) = \inf \, \{s>0 \, , \, \beta_{WLS}(s) \,
\log(\varepsilon/s) \, \leq \, 2t\} \, .$

Hence if a weak log-Sobolev inequality holds, $$\parallel P_t^*h - 1\parallel_{\L^1(\mu)} \leq
\frac{ 4 \, \int h \, \phi(h) d\mu}{(\phi \circ \tilde{\phi}^{-1})(\sqrt 2 \, \left(\int h \,
\phi(h) d\mu\right) /(\frac 1e + \varepsilon) \, \sqrt{\xi_{WLS}(\varepsilon,t)})}  \, ,$$ where
$\tilde{\phi}(u)= \sqrt{u} \phi(u)$.
\end{theorem}

The proof is analogue to the variance case. But it is shown in \cite{CGG} that:
\begin{itemize}
\item if the Poincar\'e inequality does not hold, but a weak Poincar\'e inequality holds, a weak
log-Sobolev inequality also holds (see \cite{CGG} Proposition 3.1 for the exact relationship
between $\beta_{WP}$ and $\beta_{WLS}$) but yields a worse result for the decay in total variation
distance, i.e. in this situation Theorem \ref{thmweaklogsob} is not as good as Theorem
\ref{thmweakpoincare}, \item if a Poincar\'e inequality holds, one can reinforce the weak
log-Sobolev inequality into a restricted log-Sobolev inequality.
\end{itemize}
We shall thus describe this reinforcement.

\begin{theorem}\label{thmrestricted}
Assume that $\mu$ satisfies a Poincar\'e inequality with constant $C_P$ and a weak logarithmic
Sobolev inequality with function $\beta_{WLS}$, and define $\gamma_{WLS}(u)=\beta_{WLS}(u)/u$.
Then for all $t>0$ and all bounded density of probability $h$, it holds $$\Ent_\mu(P_t^* h) \,
\leq \, e^{ - \, t / 2 \gamma_{WLS}^{-1}(\sqrt{3 C_P} \, \parallel h \parallel_\infty)} \,
\Ent_\mu(h) \, .$$ Hence if $\int h \phi(h) d\mu < + \infty$,
\begin{itemize}
\item if $\phi(u) \geq c \, \log(u)$ at infinity for some $c>0$, there exists a constant $c(\phi)$
such that
$$\parallel P_t^*h - 1\parallel_{\L^1(\mu)} \leq \frac{ c(\phi) \, \int h \, \phi(h) d\mu}{\phi \circ
\zeta_{WLS}^{-1}(t)} \, ,$$ where $\zeta_{WLS}(u)=2 \, \log(\phi(u)) \, \gamma^{-1}_{WLS}(\sqrt{3
C_P} u)$, \item if $\phi(u) \leq c \log(u)$ at infinity for all $c>0$, there exists a constant
$c(\phi)$ such that
$$\parallel P_t^*h - 1\parallel_{\L^1(\mu)} \leq \frac{ c(\phi) \, \left( 1 +\int h \, \phi(h) d\mu\right)}
{\phi \circ \theta_{WLS}^{-1}(t)} \, ,$$ where $\theta_{WLS}(u)= 2 \, \log(\phi(u) \log(u)) \,
\gamma^{-1}_{WLS}(\sqrt{3 C_P} u)$.
\end{itemize}
\end{theorem}
\begin{proof}
The first result is \cite{CGG} Proposition 4.2. We just here give the explicit expression of
$\gamma_{WLS}$. Using \eqref{eqTV} and this result give the result if we add two remarks : in the
first case we may find $C_\phi$ such that $\Ent_\mu(h\wedge K) \leq C_\phi \left(\int h \phi(h)
d\mu\right)$ for all $K > e$, so that the result follows with $c_\phi=2+ C_\phi$; in the second
case we use $\Ent_\mu(h\wedge K) \leq \log(K)$.
\end{proof}
\medskip

\subsection{Examples}

In the previous subsections, we introduce a bench of inequalities, Poincar\'e inequality or its
weak version and logarithmic Sobolev inequality and also its weak version, for which necessary and
sufficient conditions exist in dimension 1, and for which sufficient conditions are known in the
multidimensional case. Results in dimension 1 relies mainly on explicit translation of capacity
measure criterion established in \cite{BR03,BCR2,BCR1,CGG}, and we refer to their works for
further discussion. However, capacity measure conditions are (up to the knowledge of the authors)
of no use in the multidimensional setting. Let us consider the following (simplified) case: assume
that $d\mu=e^{-2V}dx$ for some regular $V$. A sufficient well known condition for a Poincar\'e
inequality to hold (see \cite{logsob} for example) is that there exists $c$ such that
$$|\nabla V|^2-\Delta V \geq c > 0$$
for large $x$'s. The associated generator is $L={1\over 2}\Delta-\nabla V.\nabla$. For general
reversible diffusion the following (nearly sufficient for exponential decay) drift condition (see
\cite{BCG} Th.2.1 or for explicit expressions of constant Th. 3.6): $\exists u\ge1$, $\alpha,b>0$
and a set $C$ such that
$$Lu(x)\le-\alpha u(x)+b1_C(x)$$
which are easy to deal with conditions which moreover extend to the weak Poincar\'e setting
\cite[Th.3.10 and Cor. 3.12]{BCG}: $\exists u\ge1$, $\alpha,b>0$, a positive function $\phi$ and a
set $C$ such that
$$Lu(x)\le -\phi(u(x))+b1_C(x).$$
As a more precise example, consider the diffusion process
$$dX_t=b(X_t)dt+\sigma(X_t)dW_t$$
where the diffusion matrix $\sigma$ has bounded smooth entries and is uniformly elliptic and assume
$$ \exists 0<p<1, M,r>0 \mbox{ such that }\forall |x|>M,\,\, x.b(x)\le -r|x|^{1-p} \, .$$
Then the invariant measure satisfies a weak Poincar\'e inequality with
$\beta_W(s)=d_p\log(2/s)^{2p/(1+p)}$ and Th. \ref{thmweakpoincare} implies for $0<q<1$
$$\|\,P^*_th-1\,\|_{\L^1(\mu)}\le C_{p,q}\left(\int h^{1+q}d\mu\right)^{1/(1+q)}e^{-D_{p,q}t^{1-p\over1+p}}.$$

There are also well known conditions for logarithmic Sobolev inequalities. Among them the most
popular is the Bakry-Emery condition: assume that $V(x)=v(x)+w(x)$ where $w$ is bounded and $v$
satisfies $Hess(v)\ge\rho Id$ for some positive $\rho$ then a logarithmic Sobolev inequality holds
with constant $e^{osc(w)}/\rho$. One may also cite Wang \cite{Wbook} and Cattiaux \cite{cat5} for
conditions in the lower bounded (possibly negative) curvature case plus integrability assumptions
or on drift like conditions. Both are however non quantitative and are thus not interesting for
our study. Concerning weak logarithmic Sobolev inequalities, in the regime between Poincar\'e and
logarithmic Sobolev inequalities, the only multidimensional conditions known can be obtained
through a $F$-Sobolev inequality we shall describe further in the next section.
%Note that Wang and
%Wu and the authors \cite{CGWW} give drift conditions for Super-Poincar\'e inequality (and thus
%logarithmic Sobolev inequality, or $F$-Sobolev inequality) which will prove effective for this
%study but this work is still in progress.

%%%%%%%%%%%%%%%%%%%%%%%%%%%%%%%%%%%%%
%%%%%%%%%%%%%%%%%%%%%%%%%%%%%%%%%%%%%
%%%%%%%%%%%%%%%%%%%%%%%%%%%%%%%%%%%%%

\section{\bf Some general results on $\Ipsi$-inequalities.}\label{seccapacity}

In this section we shall give some general results on $\Ipsi$-inequalities, i.e there exists
$C_\psi>0$ such that for all nice functions $h$
$$\int \psi(h)d\mu\le C_\psi\int\psi''(h)\Gamma(h)d\mu.$$
\smallskip

First, we use the usual way to derive Poincar\'e inequality from a logarithmic Sobolev inequality
i.e. we write $h=1+\varepsilon g$ for some bounded $g$ such that $\int g d\mu =0$. For
$\varepsilon$ going to 0 (note that $h$ is non-negative for $\varepsilon$ small enough), we see
that if $\psi$ satisfies $\psi''(1)>0$, an $\Ipsi$-inequality implies a Poincar\'e inequality
$$\Var_\mu(g) \, \leq \, 2 \, C_\psi \, \int \, \Gamma(g) \, d\mu \, ,$$ i.e. with a Poincar\'e
constant $C_P=2 C_\psi$.
\smallskip

Next, for our purpose, what is important is to control some moment of $h$. Hence what really
matters is the asymptotic behavior of $\psi$. In particular if $\eta$ is a function which is
convex at infinity (i.e. $\eta''(u)>0$ for $u\geq b$) and such that $\eta(u)/u$ goes to infinity
at infinity, we may build some ad-hoc $\psi$ as follows.

For $a>2\wedge b$, we define
\begin{eqnarray}\label{eqdefpsip}
\psi''(u) & = & \frac{\eta''(u)}{\eta''(a)} \quad \textrm{ if }
u\geq a \quad , \quad \psi''(u)  =  1 \quad \textrm{ otherwise,} \\
\psi'(u) & = & \int_{\frac 12}^u \, \psi''(v) \, dv \quad \textrm{ and } \quad \psi(u) \, = \,
\int_1^u \, \psi'(v) \, dv \, .\nonumber
\end{eqnarray}
It is easily shown that $\psi(u)= \frac 12 \, (u^2-u)$ for $u\leq a$, while one can find some
constants $\beta$ and $\gamma$ such that $\psi(u) = (\eta(u)/\eta''(a)) + \beta \, u + \gamma$ for
$u\geq a$, so that there is a constant $c$ such that $\psi(u) \leq c \, \eta(u)$ for $u\geq a$
(recall that $\eta(u)/u$ goes to infinity at infinity) and $\psi(u) \geq \frac{1}{2 \, \eta''(a)}
\, \eta(u)$ for $u$ large enough.

The choice of $\frac 12$ in the definition of the derivative, ensures that $\psi$ is non-positive
for $u\leq 1$.
\smallskip

The function $\psi$ fulfills the assumptions in Lemma \ref{lempinsker}. Of particular interest
will be the associated inequality \eqref{eqineqpsi} which, as we already remarked implies a
Poincar\'e inequality with $C_P=2 C_\psi$. We may look at sufficient conditions for
an $\Ipsi$-inequality to be satisfied.
\smallskip

\subsection{A capacity measure condition for an $\Ipsi$-inequality}

Let us first reduce the study of an $\Ipsi$-inequality to the large value case via the  use of
Poincar\'e inequality. Indeed, as previously pointed out, it is a natural assumption to suppose
that $\mu$ satisfies some Poincar\'e inequality with constant $C_P$. To prove that $\mu$ satisfies
\eqref{eqineqpsi} it is enough to find a constant $C$ such that
\begin{eqnarray*}
\int_{h\leq a} \, (h^2 - h) \, d\mu & \leq & C \, \int_{h\leq a} \, \Gamma(h) \, d\mu \,
\quad \textrm{ and }\\
\int_{h> a} \,  \eta(h) \, d\mu & \leq & C \, \left( \int_{h\leq a} \, \Gamma(h) \, d\mu +
\int_{h>a} \, \eta''(h) \, \Gamma(h) \, d\mu\right) \, .
\end{eqnarray*}
Indeed, \underline{up to the constants}, the sum of the left hand sides is greater than  $\int \,
\psi(h) \, d\mu$, while the sum of the right hand sides is smaller than $\int \, \psi''(h) \,
\Gamma(h) \, d\mu$.
\smallskip

For the first inequality, let $h$ be a nice density of probability ($h$ belongs to the domain
$\mathbb D(\Gamma)$ of the Dirichlet form $\mathcal E(h) = \int \Gamma(h) d\mu$). Remember that
$\int (h\wedge a) d\mu \leq 1$. Hence
\begin{eqnarray*}
\int_{h\leq a} \, (h^2 - h) \, d\mu & \leq & \int \, ((h\wedge a)^2 - (h\wedge a)) \, d\mu \\ &
\leq & \int \, (h\wedge a)^2 \, d\mu \, - \, \left(\int \, h\wedge a \, d\mu\right)^2 \\ & \leq &
C_P \, \int \, \Gamma(h\wedge a) \, d\mu \, = \, C_P \, \int_{h\leq a} \, \Gamma(h) \, d\mu
\end{eqnarray*}
applying the Poincar\'e inequality with $h\wedge a$ which belongs to $\mathbb D(\Gamma)$. For the
latter equality we use the second part of \eqref{eqchainrule} for a sequence $\Psi_n$
approximating $u \mapsto u\wedge a$ and use Lebesgue bounded convergence theorem.
\smallskip

To manage the remaining term, we introduce some capacity-measure condition, whose origin can be
traced back to Mazja \cite{Maz}. Following \cite{BR03,BCR1}, for $A \subset \Omega$ with
$\mu(\Omega) \leq 1/2$, we define
$$Cap_\mu(A,\Omega) \, := \, \inf \, \{ \, \int \, \Gamma(f) \, d\mu \, \, ; \, \, \BBone_A \, \leq f \,
\leq \, \BBone_\Omega \, \} \, ,$$ where the infimum is taken over all functions in the domain of the
Dirichlet form. By convention this infimum is $+\infty$ if the set of corresponding functions is
empty.

If $\mu(A) < 1/2$ we define $$Cap_\mu(A) \, := \, \inf \, \{ \, Cap_\mu(A,\Omega) \, ; \, A
\subset \Omega \, , \, \mu(\Omega) \leq 1/2 \, \} \, .$$ A capacity measure condition is usually
stated as the existence of some function $\gamma$ such that $\gamma(\mu(A))\le C Cap_\mu(A)$. Such
an inequality, and depending on the form of $\gamma$, is (qualitatively) equivalent to nearly all
usual functional inequalities: (weak) Poincar\'e inequality, (weak) logarithmic Sobolev
inequality, $F$-Sobolev inequality or generalized Beckner inequality. It is then a precious tool
to compare those inequalities, translating then properties of one to the other or using known
conditions for one to the other. It has moreover the good taste to be explicit in dimension 1. It
is then natural to look at some capacity-measure condition for an $\Ipsi$-inequality. Our first
result (similar to Theorem 20 in \cite{BCR1}) is the following

\begin{theorem}\label{thmcapacity}
Assume that $\mu$ satisfies a Poincar\'e inequality with constant $C_P$. Suppose
\begin{itemize}
\item[$\mathbf (H_\eta)$]: let $\eta$ be a $C^2$
non-negative function defined on $\R^+$ such that
\begin{itemize}
\item $\lim_{u \to +\infty}\eta(u)/u = +\infty$ , \item there exists $b>0$ such that $\eta''(u)>0$
for $u>b$, \item $\eta$ is non-decreasing on $[b, +\infty)$ and $\eta''$ is non-increasing on
$[b,+\infty)$ .
\end{itemize}
\item[$\mathbf (H_F)$]: there exist $\rho>1$ and a non-decreasing function $F$ such that
\begin{itemize}
\item for all $A$ with $0 < \mu(A)< 1/2$, $\mu(A) \, F(1/\mu(A)) \, \leq \, Cap_\mu(A)$, \item
there exists a constant $C_{cap}$ such that for all $u>a$, $$\frac{\eta(\rho \, u)}{u^2 \,
\eta''(u) \, F(u)} \, \leq \, C_{cap} \, .$$
\end{itemize}
\end{itemize}
Then $\mu$ satisfies an $\Ipsi$-inequality for $\psi$ defined in \eqref{eqdefpsip},
hence $I_\psi(t,h) \leq e^{-t/2 C_\psi} \, I_\psi(h)$. In particular there exist constants
$M_\eta$ and $C_\eta$ such that
$$\parallel P^*_t(h) \mu - \mu
\parallel_{TV} \, \leq \, M_\eta \, e^{- \, t/ 4 \, C_\eta} \, \left(1 + \int \eta(h) d\mu\right) \, .$$
\end{theorem}
\begin{proof}
According to the previous discussion, it remains to control $\int_{h> a} \,  \eta(h) \, d\mu $.
Define $\Omega =\{h>a\}$. By the Markov inequality $\mu(\Omega) \leq 1/a \leq 1/2$ since $a>2$.

For $k \geq 0$, define $\Omega_k =\{h>a \rho^k \}$ for $\rho>1$ previously defined. Again
$\mu(\Omega_k) \leq 1/(a \, \rho^k)$ and
$$\int_{h> a} \,  \eta(h) \, d\mu \, \leq \, \sum_{k \geq 0} \, \int_{\Omega_k \backslash
\Omega_{k+1}} \, \eta(h) \, d\mu \, \leq \, \sum_{k \geq 0} \, \eta(a \, \rho^{k+1}) \,
\mu(\Omega_k) \, ,$$ since $\eta$ is non-decreasing on $[a,+\infty)$. But thanks to our
hypothesis, $$\mu(\Omega_k) \leq \frac{Cap_\mu(\Omega_k)}{F(1/\mu(\Omega_k))} \leq
\frac{Cap_\mu(\Omega_k)}{F(a \, \rho^k)} \, ,$$ since $F$ is non-decreasing, provided
$\mu(\Omega_k) \neq 0$. Since $\Omega_k \supseteq \Omega_{k+1}$ the previous sum has thus to be
taken for $k < k_0$ where $k_0$ is the first integer such that $\mu(\Omega_{k_0}) = 0$ if such an
integer exists. So from now on we assume that $\mu(\Omega_k) \neq 0$.

Consider now, for $k\geq 1$ the function $$f_k := \min \, \left(1 \, , \, \left(\frac{ h \, -
 \, a \, \rho^{k-1}}{a \, \rho^{k} \, - \, a \, \rho^{k-1}}\right)_+ \right)
 \, .$$  Since $\mu(\Omega_{k-1}) < 1/2$ and $f_k$ vanishes
 on $\Omega_{k-1}^c$, $f_k$ vanishes with probability at least $1/2$. Hence $$Cap_\mu(\Omega_k) \,
 \leq \, \int \, \Gamma(f_k) \, d\mu \, \leq \, \frac{\int_{\Omega_{k-1} \backslash \Omega_k} \,
 \Gamma(h) \, d\mu}{a^2 \, \rho^{2(k-1)} \, (\rho - 1)^2} \, \leq \,
 \frac{\int_{\Omega_{k-1} \backslash \Omega_k} \, \eta''(h) \,
 \Gamma(h) \, d\mu}{a^2 \, \rho^{2(k-1)} \, (\rho - 1)^2 \, \eta''(a \, \rho^{k})} \, ,$$ since
 $\eta''$ is non-increasing.

Summing up all these estimates (for $k\geq 1$ remember) we obtain
\begin{eqnarray}\label{eqinterk}
\int_{h>\rho \, a} \, \eta(h) \, d\mu & \leq & \sum_{k \geq 1}  \left(\frac{\eta(a \,
\rho^{k+1})}{a^2 \, \rho^{2(k-1)} \, (\rho - 1)^2  \eta''(a \, \rho^{k})  F(a \, \rho^k)}\right)
\, \int_{\Omega_{k-1} \backslash \Omega_k}  \eta''(h) \,
 \Gamma(h) \, d\mu \nonumber \\ & \leq & \frac{\rho^2 \, C_{cap}}{(\rho - 1)^2} \, \int_{h> a} \, \eta''(h) \,
 \Gamma(h) \, d\mu \, ,
 \end{eqnarray}
according to our hypothesis.

It remains to control $\int_{a< h \leq \rho \, a} \, \eta(h) \, d\mu$. But on $\{a< h \leq \rho \,
a\}$, $\eta(h) \leq c \, (h^2-h)$ for some $c>0$, and as before $$\int_{h<a \rho} (h^2-h) d\mu \,
\leq \, \int_{h<a \rho} \, \Gamma(h) \, d\mu \leq C' \, \left( \int_{h\leq a} \, \Gamma(h) \, d\mu
+ \int_{h>a} \, \eta''(h) \, \Gamma(h) \, d\mu\right)$$ for some $C'$ since $\eta''$ is bounded
from below on $[a, a \, \rho]$. The proof is completed.
\end{proof}

\begin{remark}\label{remgeneralites} \textbf{Remarks and examples.}
\begin{enumerate}
\item[(1)] \quad If $\eta(u)=u^2$ we may choose $F(u)=c$ for all $u$ and conversely. The
capacity-measure inequality $\mu(A) \leq (1/c) \, Cap_\mu(A)$ is known to be equivalent (up to the
constants) to the Poincar\'e inequality. We thus recover (see below for more precise results) the
usual $\L^2$ theory. Note that as we suppose $F$ to be non decreasing, so that $\mathbf(H_F)$
already implies a Poincar\'e inequality, but with no precision on the constant.

Similarly if $\eta(u)= u \log(u)$ we may choose $F(u)= C \log(c u)$ for some well chosen $c$, $C$
and conversely. Again the capacity-measure inequality $\mu(A) \, \log(c/\mu(A)) \, \leq (1/C) \,
Cap_\mu(A)$ is known to be equivalent (up to the constants) to the logarithmic Sobolev inequality,
and we recover the usual entropic theory.

\item[(2)] \quad Since we know now what hypotheses on $\eta$ are required we may follow more
accurately the constants. Indeed since $\eta''$ is non-increasing, $\psi'' \leq 1$ for $u>a$. It
is thus not difficult to check that $\psi(u) \leq (1+(\rho - 1)^2)(u^2 - u)$ on $[a, \rho \, a]$
(using $a>2$). So it easily follows that
$$C_\eta \, = \, C_\psi \, \leq \, \max \, \left(\frac{\eta''(a) \, (1+(\rho - 1)^2) \, C_P}{\eta''(\rho \,
a)} \, , \, \frac{\rho^2 \, C_{cap}}{(\rho - 1)^2}\right) \, .$$

\item[(3)] \quad Now choose $\eta(u)=u^{p}$ for some $2 \geq p>1$ (recall that $\eta''$ is
non-increasing). Again the best choice of $F$ is a constant. More precisely choose $F=3 C_P$. It
is known (see the lower bound of Theorem 14 in \cite{BCR1}) that $\mu(A) \leq F \, Cap_\mu(A)$.
Hence we have $C_{cap}= (a^2 \, \rho^{2p} / p(p-1))$. Then a rough estimate is
$$C_\eta \, \leq \, C_P \, \max \, \left(\rho^{2-p} \, (1+(\rho-1)^2) \, , \, \frac{\rho^{2+p}}{p
\, (p-1) \, (\rho-1)^2}\right) \, .$$ Hence we obtain $$\parallel P^*_t(h) \mu - \mu
\parallel_{TV} \, \leq \, M_\eta \, e^{- \, c_p \, t/  C_P} \,
 \left(1 + \int \eta(h) d\mu\right) \, ,$$ for some constant $c_p$.

If $p$ is close to one, it is easily seen that $c_p \geq (p-1) c$ for some universal constant $c$.
So, Theorem \ref{thmcapacity} explains why the results in Example \ref{examplepower} (again with a
bad constant $c$ in the previous exponential) are not so surprising.
\smallskip

A similar study is possible for $\eta(u) = u \log_+^\beta(u)$ for $\beta>0$. In this case indeed,
it is easily seen that one may choose $$F(u) \, =  \, \log(u) \quad \textrm{ and } \quad C_{cap}=
C(a,\rho) \, \frac{1\wedge 2^{\beta - 1}}{\beta} \, , $$ at least for $u$ small enough. Such a
capacity-measure is known to be equivalent to a logarithmic Sobolev inequality, and as before for
$0<\beta \leq 1$ we recover the results in Example \ref{examplepentropy} (with the linear
dependence in $\beta$ for $\beta$ close to 0).

Interesting here is also the case $\beta > 1$. Indeed one could expect that the exponential decay
of such a $\beta$-entropy would require a weaker inequality than the log-Sobolev inequality. It
seems that this is not the case, even if, as we said, we cannot claim that the $F$ obtained in
Theorem \ref{thmcapacity} furnishes the best capacity-measure inequality. \item[(4)] One may be
surprised of the intervention of a new function $F$, in $\mathbf(H_F)$, rather than  an usual
capacity-measure condition. In fact, it enables us to relax the assumptions on $\eta$. In
particular, if there exists $a$ and $\rho>1$ such that for $u>a$, $\eta(\rho u)/(u^2\rho''(u))$ is
non decreasing then instead of $\mathbf(H_F)$ one may use the capacity-measure condition: there
exists $C_{c}$ such that
$${\eta''(1/\mu(A))\over\mu(A)\,\eta(\rho/\mu(A))}\le C_{c} \, Cap_\mu(A) .$$\hfill $\diamondsuit$
\end{enumerate}
\end{remark}

\begin{remark}\label{remFSob}
Theorem \ref{thmcapacity} allows to cover the class of $F$-Sobolev inequalities. Indeed
combining the results in section 5 of \cite{BCR1} and Lemma 17 in \cite{BCR3}, if $\mu$ satisfies
a Poincar\'e inequality and the $F$-Sobolev inequality
\begin{equation}\label{eqdeffsob}
\int \, f^2 \, F\left(\frac{f^2}{\int f^2 d\mu}\right) \, d\mu \, \leq \, C \, \int \, \Gamma(f)
\, d\mu
\end{equation}
for all nice $f$, then the capacity-measure inequality in Theorem \ref{thmcapacity} is satisfied,
provided $u \mapsto F(u)/u$ is non-increasing  and $F(\lambda u) \leq (\lambda/4) \, F(u)$ for
some $\lambda>4$ and all $u$ large enough (Theorem 22 and Remark 23 in \cite{BCR1}).

Conversely, Theorem 20 in \cite{BCR1} tells us that the capacity-measure inequality in Theorem
\ref{thmcapacity} implies the $\tilde{F}$-Sobolev inequality with
$\tilde{F}(u)=\left(F(u/\rho)-F(2)\right)_+$ for $\rho>1$. With the previous hypotheses on $F$,
and up to the constants, we may replace $\tilde F$ by $F_+$.

For instance if, for $1 \leq \alpha \leq 2$, we choose $F(u)= \log^{2(1-\frac 1\alpha)}(1+u) -
\log^{2(1-\frac 1\alpha)}(2)$ the Boltzmann measure $\mu(dx)= (1/Z) e^{-2 U(x)} dx$ with
$U(x)=|x|^{\alpha}$ for large $x$, satisfies a $F$-Sobolev inequality (see \cite{BCR1} section 7).
An elementary calculation shows that we can choose $$\eta(u)= u \, \log^{2(1-\frac 1\alpha)}(u) \,
\, e^{\log^{(2/\alpha) - 1}(u)} \, ,$$ for large $u$. We thus get an interpolation result between
Poincar\'e and Gross inequalities. \hfill $\diamondsuit$
\end{remark}

\subsection{Links between $\Ipsi$-inequalities and $F$-Sobolev inequalities}

In view of the previous remark it is natural to relate an $\Ipsi$-inequality and $F$-Sobolev
inequalities. To this end define
\begin{equation}\label{eqdefH}
H(u) = \int_0^u \, \sqrt{\psi''(s)} \, ds
\end{equation}
which is a continuous increasing function, whose inverse function is denoted by $H^{-1}$. We
assume that $H(u) \to + \infty$ as $u \to + \infty$ so that $H^{-1}$ is everywhere defined on
$\R^+$. Remark that the derivative of $\psi \circ H^{-1}$ is equal to $(\psi'/\sqrt{\psi''})\circ
H^{-1}$, so is non-decreasing if $\psi''$ is non-increasing, that is $\psi \circ H^{-1}$ is a
convex function.

For $f \geq 0$, denote by
\begin{equation}\label{eqdefN}
N(f) = \inf \{ \lambda>0 \, ; \, \int H^{-1}(f/\lambda) d\mu \leq 1 \} \, .
\end{equation}
Then an easy change of variables shows that an $\Ipsi$-inequality is equivalent to
\begin{equation}\label{eqpsirondasobolev}
N^2(f) \, \int \, \psi\left(H^{-1}\left(\frac{f}{N(f)}\right)\right) \, d\mu \, \leq \, C_\psi \,
\int \, \Gamma(f) \, d\mu \, ,
\end{equation}
for all nice $f \geq 0$. \eqref{eqpsirondasobolev} looks like a $F$-Sobolev inequality except that
the normalization is not the $\L^2$ norm but $N$. As before, up to the constants, both coincide if
$F=\log$ explaining why entropy is particularly well suited.

We see that \eqref{eqpsirondasobolev} is exactly
\begin{equation}\label{eq??}
 \int
f^2 \, F(f^2/N^2(f)) \, d\mu \leq C_\psi \int \, \Gamma(f) \, d\mu \quad \textrm{ for } \quad
F(u)=(\psi \circ H^{-1})(\sqrt u)/u \, .
\end{equation}
We can thus get immediate comparison results, assuming that $F$ is non-decreasing (we will see in
the proof of the next Theorem that one can always modify \eqref{eq??} for this property to hold).
Indeed we have two interesting cases (at least for large $u$ and up to constants):
\begin{eqnarray}\label{eqdiscuss}
\textrm{either } \quad H(u)\geq \sqrt u & \Leftrightarrow u^2 \geq H^{-1}(u) \Leftrightarrow \int
f^2 d\mu \geq N^2(f) \\  \textrm{ or } \quad H(u)\leq \sqrt u & \Leftrightarrow u^2 \leq H^{-1}(u)
\Leftrightarrow \int f^2 d\mu \leq N^2(f)
\end{eqnarray}
since $H$ and $H^{-1}$ are non-decreasing. In the first case, \eqref{eq??} implies the $F$-Sobolev
inequality \eqref{eqdeffsob} while in the second case the $F$-Sobolev inequality implies
\eqref{eq??}. Note that once again the limiting case $H(u)=\sqrt u$ corresponds to log-Sobolev.

The first case gives some converse to Theorem \ref{thmcapacity}. Note that $\psi(v) = H^2(v) \,
F(H^2(v)) \geq H^2(v) \, F(v)$ since $F$ is non-decreasing, hence we get a $F$-Sobolev inequality
for some $F$ such that $F(v)\leq \psi(v)/H^2(v)$. With some additional (but reasonable)
assumptions we can improve this result. Indeed

\begin{theorem}\label{thmconvcapacity}
Let $\eta$ and $\psi$ be as in Theorem \ref{thmcapacity}, and $H$ defined in \eqref{eqdefH}. We
assume that $H(+\infty)=+\infty$.  Assume in addition that for $u$ large
\begin{itemize}
\item $u \mapsto \bar F(u) = (\psi/H^2)(u)$ is non-decreasing and satisfies $\bar F(\lambda u)\leq
\lambda \bar F(u)/4$, for some $\lambda > 4$, \item $u \mapsto \bar F(u)/u$ is non-increasing.
\end{itemize}
If $\mu$ satisfies a Poincar\'e inequality with constant $C_P$ and an $\Ipsi$-inequality for some
$C_\psi$, then for $\mu(A)$ small enough, the capacity-measure inequality $$\mu(A) \, \bar
F(1/\mu(A)) \, \leq \, D \, Cap_\mu(A) $$ is satisfied for some $D>0$. Accordingly (see Remark
\ref{remFSob}) $\mu$ satisfies the $\bar F_+$-Sobolev inequality (with some constant $D_F$).

Conversely if $\mu$ satisfies a Poincar\'e inequality with constant $C_P$ and the $\bar F$-Sobolev
inequality, and if $H(u) \geq \sqrt u$ for large $u$, an $\Ipsi$-inequality is satisfied for some
$C_\psi$.
\end{theorem}
\begin{proof}
Note that $\lim_{u \to +\infty} \bar F(u)/u$ exists by monotonicity. Denote it by $m$. We have
$\bar F(u)/4 u \geq \bar F(\lambda u)/(\lambda u)$ so that letting $u$ go to infinity we get $m/4
\geq m$ hence $m=0$. In particular the capacity-measure inequality when $\mu(A)=0$ reduces to
$Cap_\mu(A) \geq 0$ which is of course satisfied. We shall thus assume now that $\mu(A)>0$.

First we write \eqref{eq??} in the form
\begin{equation}\label{eq???}
 \int
f^2 \, F(f/N(f)) \, d\mu \leq C_\psi \int \, \Gamma(f) \, d\mu \quad \textrm{ for } \quad
F(u)=(\psi \circ H^{-1})(u)/u^2 \, .
\end{equation}

The first part of the proof is mimicking the proof of Lemma 17 in \cite{BCR3}.  Note that the
derivative of $F$ (defined in \eqref{eq???}) is given by $$u \mapsto \frac{u \, \psi'(H^{-1}(u))
\, - \, 2 \, (\psi \, \sqrt{\psi''})(H^{-1}(u))}{u^3 \, \sqrt{\psi''(H^{-1}(u))}}$$ which is
non-negative for $u$ large enough since $u \mapsto (\psi/H^2)(u)$ is non-decreasing.

Choose some $\rho >1$ large enough, such that $F(2 \rho) \geq 0$ and define $\tilde{F}(u)=
F(u)-F(2 \rho)$ which is thus non-negative for and non-decreasing on $[2\rho, +\infty)$
if $\rho$ is large enough. We thus have
\begin{equation*}
\int \, f^2 \, \tilde{F}_+\left(\frac{f}{N(f)}\right) \, d\mu \, \leq \, C_\psi \, \int \,
\Gamma(f) \, d\mu \, + \, M \, \int f^2 d\mu \, ,
\end{equation*}
with $M = \sup_{0 \leq u \leq 2\rho} |F(u)|$.

Now we can follow \cite{BCR3} with some slight modifications. We give the details for the sake of
completeness. Let $\chi$ defined on $\mathbb R^+$ as follows : $\chi(u)=0$ if $u\leq  2$,
$\chi(u)= u$ if $u\geq 2\rho$ and $\chi(u)=2\rho \, (u- 2)/(2\rho) -  2)$ if $ 2 \leq u \leq 2
\rho$. Since $\chi(f)\leq f$, $N(\chi(f)) \leq N(f)$ so that since $\tilde F_+$ is non-decreasing,
\begin{eqnarray*}
\int f^2 \tilde F_+(f/N(f)) d\mu & = & \int \chi^2(f) \tilde F_+(\chi(f)/N(f)) d\mu \\ & \leq &
\int \chi^2(f) \tilde F_+\left(\frac{\chi(f)}{N(\chi(f))}\right) d\mu \\ & \leq & B C_\psi \int
\Gamma(f) \, d\mu + M \int \chi^2(f) d\mu \\ & \leq & B C_\psi \int \Gamma(f) \,  d\mu + M
\int_{f^2 \geq 2\int f^2 d\mu} f^2 d\mu
\end{eqnarray*}
where $B=(\rho/\rho -1)^2$. But as shown in \cite{BCR1}, $\int_{f^2 \geq 2\int f^2 d\mu} f^2 d\mu
\leq 12 C_P \int \, \Gamma(f) \, d\mu$ so that we finally obtain the existence of $D_\psi$ such
that
\begin{equation}\label{eqpsirondasobolevplus}
\int \, f^2 \, \tilde{F}_+\left(\frac{f}{N(f)}\right) \, d\mu \, \leq \, D_\psi \, \int \,
\Gamma(f) \, d\mu \, .
\end{equation}
\smallskip

 The second part of the proof is mimicking the one of Theorem 22 in \cite{BCR1}. Let
$\mu(A) < 1/2$ and $\BBone_A \leq f \leq \BBone_\Omega$ with $\mu(\Omega) \leq 1/2$. For $k \in
\N$ we define $\Omega_k = \{ f \geq 2^k \, N(f)\}$ and $$f_k = \min \, \left((g - 2^k \, N(f))_+
\, ; \, 2^k \, N(f)\right) \, .$$ Note that $f_k$ is equal to 0 on $\Omega_k^c$ and to $2^k \,
N(f)$ on $\Omega_{k+1}$.

In addition, since $H^{-1}(0)=0$,
\begin{eqnarray*}
\int \, H^{-1}\left(\frac{f_k \, H(1/\mu(\Omega_k))}{2^k \, N(f)}\right) d\mu & = &
\int_{\Omega_k} \, H^{-1}\left(\frac{f_k \, H(1/\mu(\Omega_k))}{2^k \, N(f)}\right) d\mu \\ & \leq
& \int_{\Omega_k} \, H^{-1}\left(H(1/\mu(\Omega_k))\right) d\mu \, = \, 1
\end{eqnarray*}
so that $N(f_k) \leq 2^k N(f) / H(1/\mu(\Omega_k))$.  Therefore, applying
\eqref{eqpsirondasobolevplus} (we need here a non-negative $F$)
\begin{eqnarray*}
D_\psi \, \int \, \Gamma(f) \, d\mu & \geq & D_\psi \, \int \, \Gamma(f_k) \, d\mu \, \geq \,
\int_{\Omega_{k+1}} f_k^2 \, \tilde{F}_+\left(\frac{f_k}{N(f_k)}\right) \, d\mu \\ & \geq &
\mu(\Omega_{k+1}) \, 2^{2k} \, N^2(f) \, \tilde F_+(H(1/\mu(\Omega_k)) \, .
\end{eqnarray*}
We are thus in the situation of the proof of Theorem 22 in \cite{BCR1} replacing $\mu(g^2)$
therein by $N^2(f)$ and $F$ therein by $\tilde F_+\circ H$. We may conclude since for $u$ large,
$(\tilde F_+\circ H)(u) \geq c \, \psi(u)/H^2(u)$ and for $\mu(A)$ small enough according to
Remark 23 in \cite{BCR1}.
\smallskip

The direct part being proven let us briefly indicate how to prove the converse part. Again we may
modify $\bar F$ into a non-negative $G$ thanks to Poincar\'e inequality (this is exactly Lemma 17
in \cite{BCR3}). The properties of $G$ ensure that we may apply Theorem 22 in \cite{BCR1}, i.e.
the $G$-Sobolev inequality implies a capacity-measure inequality (with the same $G$). Next just
remark that the proof of Theorem 20 in \cite{BCR1} applies to any homogeneous inequality (i.e. we
may replace $\int f^2 d\mu$ therein by $N^2(f)$ for example). We thus get that \eqref{eq??} holds
with $G$ in place of $F$. But as we remarked $F \leq \bar F$ for large $u$, and with our
hypotheses $\bar F \leq c \, G$ at infinity. We may thus replace (changing the constants) $G$ by
$F$ for large $u$, small values of $u$ can be controlled again (if necessary) by using Poincar\'e
inequality.
\end{proof}

\begin{remark}\label{remlemeilleur}
At least if $H(u) \geq \sqrt u$ (up to a constant actually), we have two results saying that some
$F$-Sobolev inequality implies an $\Ipsi$-inequality: the first one with $F(u) \geq C \, \eta(\rho
u)/u^2 \, \eta''(u)$ at infinity, the second one with $\bar F(u)= \eta(u)/H^2(u)$. It seems not
easy to compare them in full generality. However one can use some asymptotic estimates.

First recall that $\psi''$ (hence $\sqrt{\psi''}:=g$) is supposed to be non-increasing at
infinity. Since we have assumed that $H(+\infty)=+\infty$ it implies that $g'(u)/g(u) = (1/2)
(\psi'''(u)/\psi''(u)) \, \geq \, - 1/u$ near infinity. Now write the elementary $$\int_m^u \,
(g(s) + s g'(s)) \, ds \, = \, u g(u) \, - \, mg(m) \, .$$ It immediately follows that
\begin{equation}\label{eqcomparaisondieudonne}
\textrm{if } \, \frac{u \, \psi'''(u)}{\psi''(u)} \, \to \, 0 \, \textrm{ as } \, u \to +\infty \,
, \, \textrm{ then } \, H(u) \sim_{u \to +\infty} u \, \sqrt{\psi''(u)} \, ,
\end{equation}
while
\begin{equation}\label{eqcomparaisondieudonnebis}
\textrm{if } \, \liminf_{u \to +\infty} \, \frac{u \, \psi'''(u)}{2 \, \psi''(u)} \, =\, - \, d \,
 , \, \textrm{ for some $d<1$, then } \, H(u)  \leq _{u \to +\infty} \, \frac{1}{1-d} \, u \, \sqrt{\psi''(u)}
 \, .
\end{equation}
Hence we always get that
\begin{equation}\label{eqalleluia}
\bar F(u) \, \geq c \, \frac{\psi(u)}{u^2 \, \psi''(u)} \, ,
\end{equation}
that is in general the same condition in both Theorems. This is very satisfactory but of course we
have made additional assumptions on $\bar F$ in Theorem  \ref{thmconvcapacity}. \hfill
$\diamondsuit$
\end{remark}

One of the very interesting feature of $F$-Sobolev inequality is that they are linked to
contraction properties for the semi-group. We now recall these general results taken from
\cite{Wbook}.
\smallskip

According to Wang's beautiful results (\cite{Wbook} chapter 3.3), a $F$-Sobolev inequality is
equivalent to a super-Poincar\'e inequality, i.e. for all nice $f$ and all $s\geq 1$,
\begin{equation}\label{eqsuperpoinc}
\int f^2 d\mu \leq \beta_{SP}(s) \, \int \Gamma(f) d\mu + s \left(\int |f| d\mu\right)^2 \, .
\end{equation}
If the $F$-Sobolev inequality holds, \eqref{eqsuperpoinc} holds with $\beta_{SP}(s)= c/F(s)$ for
$s$ large enough (\cite{Wbook} Theorem 3.3.1). For a somewhat intricate converse see \cite{Wbook}
Theorem 3.3.3.

Assume that the $F$-Sobolev inequality holds. The associated super-Poincar\'e inequality implies
some boundedness for the associated semi-group. Of particular interest here are Theorem 3.3.13 (2)
and Theorem 3.3.14 in \cite{Wbook}. The first one tells us that $P_t$ is super-bounded (i.e. is
bounded from $\L^2(\mu)$ in $\L^p(\mu)$ for all $p>2$ and all $t>0$) as soon as $F(u)/\log(u) \to
+\infty$ as $u \to \infty$ (some converse statement is also true), while the second one tells us
that $P_t$ is ultracontractive (or ultrabounded in Wang's terminology) as soon as $$\int^{+\infty}
\, \frac{1}{u \, F(u)} \, du \, < \, +\infty \, .$$
\smallskip

Let us come back to the second situation in \eqref{eqdiscuss}. Roughly speaking this case is the
one of stronger inequalities than the log-Sobolev inequality, for which with the mild additional
previous assumptions, we know that the semi-group is ultracontractive. However we can give another
interesting example, and will continue the discussion in the next section.

\begin{example}\label{exsuperbounded}
For $F(u)= \log(u) \, \log(\log(u))$ at infinity, Wang's results show that the semi-group $P_t^*$
is super-bounded but not ultracontractive. An elementary calculation show that we can choose
$\eta(u) = u \, \log(\log(u))$ in this case. \hfill $\diamondsuit$
\end{example}
\medskip

The study of weak inequalities should be interesting. The two extreme cases, weak Poincar\'e and
weak logarithmic Sobolev inequalities have already been studied. As remarked in \cite{CGG} the
main interest of weak log-Sobolev inequalities is to describe some interpolation between
Poincar\'e and Gross (if a Poincar\'e inequality does not hold, the weak log-Sobolev inequality
furnishes worse results than the corresponding weak Poincar\'e inequality). So the potential weak
inequalities should give better results than the weak log-Sobolev inequality (recall Theorem
\ref{thmrestricted}). However, the technical intricacies are certainly too much for a potential
reader since we do not have (yet) any convincing application.

\begin{remark}\label{remarno}
Finally we may ask whether it is possible to get some exponential decay using a weaker inequality
than Poincar\'e inequality but for $\eta$'s larger than $u \mapsto u^2$ at infinity.

Assume for instance that for all density of probability $h$ bounded by $M\geq 2$ we have for some
function $\xi$ decaying to 0, $$\int |P_t^*h - 1| d\mu \leq \xi(t) \, .$$ Let $f$ be in
$\L^2(\mu)$ such that $\int f d\mu = 0$ and $\parallel f \parallel_\infty \leq 1$. Then
$h=(f+2)/2$ is a density of probability, bounded by $3/2$ hence
$$
\Var_\mu(P_t^* f) \, \leq \, 2 \, \int |P_t^*h - 1| d\mu \, \leq \, 2 \, \xi(t) \, \leq \, 2 \,
\xi(t) \, \Osc^2(f)
$$
and the previous inequality extends to all $f$ in $\L^2(\mu)$ by homogeneity.

In the symmetric case ($P_t=P_t^*$) this result implies a weak Poincar\'e inequality (see
\cite{RW01} Theorem 2.3). In particular if $\xi(t)=c \, e^{- \lambda t}$ for some $\lambda >0$ the
same Theorem shows that $\mu$ satisfies a Poincar\'e inequality. Hence in the symmetric case we
cannot obtain any exponential decay for the total variation distance even for bounded densities
without assuming that a Poincar\'e inequality is satisfied. If it is not we have to use the
results of the previous section.\hfill $\diamondsuit$
\end{remark}
\begin{remark}
An aficionado of functional inequalities may have remarked that we have not discussed usual
properties introduced when dealing with a new functional inequality like $\Ipsi$: tensorization
and concentration of measure. In fact, concentration is not at all our purpose here and in fact it
may be directly deduced from the capacity measure condition imposed in Theorem \ref{thmcapacity}
or inherited by Theorem \ref{thmconvcapacity}. Concerning tensorization, it is more relevant for
applications concerning diffusion to deal directly in multidimensional space rather than the
limiting setting of tensorization and perturbation argument. Note also that  by the equivalence
obtained via Theorem \ref{thmconvcapacity}, of an $\Ipsi$ inequality and an $F$-Sobolev
inequality, we get all the tensorization property (and concentration) via $F$-Sobolev
inequalities, see \cite{BCR1,BCR3} for details.
\end{remark}

\subsection{Further examples}

The major difference between  Theorems \ref{thmcapacity} and \ref{thmconvcapacity} is that in  the
first one we do not explicitly suppose an $F$-Sobolev inequality. Therefore we may put less
stringent assumptions on $F$, and still have an explicit condition in dimension 1: namely
$\mathbf(H_F)$ can be translated in
\begin{itemize}
\item[$\mathbf (H'_F)$]: there exist $\rho>1$ and a non-decreasing function $F$ such that
\begin{itemize}
\item let $m$ be a median of $\mu$, and denoting $\mu_c$ the density of the absolutely continuous
part of $\mu$ w.r.t. the Lebesgue measure, if
$$\sup_{x>m}\mu([x,\infty[)F(1/\mu([x,\infty[))\int_m^x\mu^{-1}_c(t)dt\,<\,\infty$$
$$\sup_{x<m}\mu(]-\infty,x])F(1/\mu(]-\infty,x]))\int_x^m\mu^{-1}_c(t)dt\,<\,\infty$$
\item
there exists a constant $C_{cap}$ such that for all $u>a$, $$\frac{\eta(\rho \, u)}{u^2 \,
\eta''(u) \, F(u)} \, \leq \, C_{cap} \, .$$
\end{itemize}
\end{itemize}
However this measure capacity condition is no more tractable in the multidimensional case whereas
we have known conditions in the multidimensional case for $F$-Sobolev inequalities. Indeed, by
\cite[Th. 21]{BCR3}, assume that $d\mu=e^{-2V}dx$ with $V$ a $C^2$ potential such that $Hess(V)\ge
R$ for some real $R$ and let $F$ be $C^1$ on $]0,\infty[$ such that
\begin{itemize}
\item $F(x)\to\infty$ as $x\to\infty$, $F(x)\le c\log_+x$, $F(xy)\le \hat c+F(x)+F(y)$ and
$xF'(x)\le \tilde c$ for some positive $c,\tilde c$ and real $\hat c$; \item the following drift
like condition is verified: $F(e^{2V})+C(LV-|\nabla V|^2)\le K$ for some positive $C$ and $K$
\end{itemize}
then $\mu$ verifies a ($F$-$B$)-Sobolev inequality for some positive $B$. Using then  Theorem
\ref{thmcapacity} via \cite[Th. 18]{BCR3} for the implied capacity-measure condition $\mathbf
(H_F)$ we get an $\Ipsi$-inequality, hence an exponential decay for the total variation distance
using Lemma \ref{lempinsker}.

Consider for example, for $1<\alpha<2$, $V(x)=|x|^\alpha +\log(1+|x|\sin^2(x))$, then $\mu$
satisfies  a Poincar\'e inequality and the previous conditions with
$F(u)=\log(1+u)^{2(1-1/\alpha)}-\log(2)^{2(1-1/\alpha)}$, so that we get for some $c_1,c_2>0$
$$\|\, P^*_th-1\,\|_{\L^1(\mu)}\le c_1\,e^{-c_2t}\left(\int h\log^{2(1-1/\alpha)}
(h)e^{\log^{(2/\alpha)-1}(h)}d\mu)\right)^{1/2}.$$
%For more refined conditions ensuring existence of $F$-Sobolev inequalities via Super-Poincar\'e
%inequalities, we refer to \cite{CGWW}.

%%%%%%%%%%%%%%%%%%%%%%%%%%%%%%%%%%%%%%%
%%%%%%%%%%%%%%%%%%%%%%%%%%%%%%%%%%%%%%%
%%%%%%%%%%%%%%%%%%%%%%%%%%%%%%%%%%%%%%%

\section{\bf Is a direct study of the total variation distance possible ?}\label{secTV}

A natural question is of course : is it possible to directly study the possible decay of the total
variation distance, instead of looking at larger quantities like the variance or the relative
entropy ? Due to the non smoothness of $u \mapsto |u-1|$ the answer is no, but one can try to
replace the total variation distance by almost equivalent quantities.
\smallskip

Before to look at such cases, we just make a remark. It is an easy consequence of the semi-group
property that, if $$B^1=\{f\in \L^1(\mu) \, ; \, \int f d\mu = 0 \, , \, \int |f| d\mu \leq 1\}$$
an uniform decay
$$\sup_{f \in B^1}\parallel P_{t_0}^* f\parallel_{\L^1(\mu)} \leq e^{-\lambda} < 1$$
for some $t_0>0$ implies an exponential decay $$\sup_{f \in B^1}\parallel P_{t}^*
f\parallel_{\L^1(\mu)} \leq C \, e^{-\lambda \, t} \, .$$ A similar result for the total variation
distance $$\parallel (P_{t_0}^*h)\mu - \mu\parallel_{TV} \leq e^{-\lambda} $$ implies an
exponential decay is not clear. Of course if we assume that for all $h$, $$\parallel
(P_{t_0}^*h)\mu - \mu\parallel_{TV} \leq e^{-\lambda} \,
\parallel h\mu - \mu\parallel_{TV}$$ the semi-group property again implies an exponential decay.

However one can suspect that such an uniform decay for the total variation distance is a very
strong result.

The Ornstein-Uhlenbeck process (on $\R$) for instance does not satisfy this property since the law
at time $t$ starting from $x$ is given by
$$P_t(x,dy) = \left(\pi (1-e^{-2t})\right)^{- 1/2} \, \exp - \, \frac{(y-xe^{-t})^2}{1 - e^{-2t}}
\, dy \, ,$$ so that if $\mu(dy)$ is the Gaussian measure with zero mean and variance $1/2$ (which
is reversible for the process), for large $t>0$, choosing $x= e^t$, we obtain $$ \sup_x \parallel
P_t(x,.)-\mu \parallel_{TV} \geq \sqrt{1/2\pi} \, \int \, |e^{- \left((y-1)^2 - y^2\right)} - 1|
\, \mu(dy) \geq c > 0 \, .$$

Of course if $P_t^*$ is ultracontractive we have an exponential decay in $\L^1$. So extending the
result of the preceding section to the linear case, should have some interest in the study of
ultracontractivity.

\subsection{\bf The linear case.}\label{seclinear}

In the preceding section we only looked at functions $\psi$ such that $\psi(u)/u \to \infty$ at
infinity. However Remark \ref{remlinear} shows that it is possible to consider cases where $\psi$
is almost linear at infinity.

Consider (at least for large $u$) $\eta(u)= u + \theta(u)$ where $\theta$ is a convex function
such that $\theta(u)/u \to 0$ as $u \to \infty$. Necessarily $\theta'(u) \leq 0$ for large $u$ and
goes to 0 as $u \to \infty$. If $\eta''$ is non-increasing, so does $\theta''$, and according to
the previous property $\theta''(u) \to 0$ at infinity.

Define $\psi$ as in \eqref{eqdefpsip}. Then for $u>a$ , $\psi(u)= \frac{- \theta'(a)}{\theta''(a)}
\, u \, + \, \nu(u)$ where $\nu(u)/u \to 0$ at infinity. In order to apply Remark \ref{remlinear}
it is thus enough to have $2 \theta'(a)+\theta''(a) < 0$ (since $\psi'(1)=1/2$).

Assuming this condition, we may extend Theorem \ref{thmcapacity} to this $\eta$. This yields $F(u)
\geq \eta(\rho \, u)/(u^2 \, \theta''(u))$ and since what is important is the behavior of $F$ near
infinity and $\eta$ is moderate, the key is the behavior of $u \mapsto 1/(u \, \theta''(u))$ when
$u$ goes to infinity. A capacity-measure inequality is interesting only if $ F(u)/u \to 0$ as $u$
goes to infinity (otherwise we already know that the semi-group is ultracontractive) so that the
only interesting cases are those for which $u^2 \, \theta''(u) \to \infty$ as $u \to \infty$.

The main question is: is it possible to build such $\theta$'s ? The simplest way to do so is to
write
\begin{equation}\label{eqtheta}
\theta(u)= - \, \int_a^u \, (1/\tau(s)) ds \quad , \quad \theta'(u)= - \, (1/\tau(u)) \quad ,
\quad \theta''(u) = (\tau'(u)/\tau^2(u))
\end{equation}
where $\tau$ is a non-negative, non-decreasing function. Fix some $F$. In our situation what we
have to do is to find some $\tau$ such that
$$\frac{\tau'(u)}{\tau^2(u)} \, = \, \frac{1}{u \, F(u)} \, .$$ Since $\theta'=-1/\tau$ goes to 0
at infinity, it implies Wang's integrability condition, hence ultracontractivity.
\smallskip

Since $c I_\psi(h) \leq \parallel h - 1\parallel_{\L^1(\mu)} \leq C \sqrt{I_\psi(h)}$ because
$\psi$ is almost linear at infinity, an exponential decay of the total variation distance
($\parallel P^*_t h - 1\parallel_{\L^1(\mu)} \leq C \, e^{-\alpha t}$) is equivalent to the
exponential decay of $I_\psi(t,h)$, but with an initial control by $\sqrt{I_\psi(h)}$. What we
just did is to show that such an exponential decay $I_\psi(t,h) \leq e^{-\alpha t} \, I_\psi(h)$
(notice that this inequality is an equality at $t=0$) cannot be obtained through a $F$-Sobolev
inequality, unless $P_t$ is ultracontractive. However Theorem \ref{thmcapacity} only furnishes one
direction : $F$-Sobolev implies uniform exponential decay. So we cannot claim, but we strongly
suspect that the uniform exponential decay of the total variation distance is actually equivalent
to ultracontractivity.
\medskip

\subsection{\bf Using the Hellinger distance.}\label{sechellinger}

Another possibility to control the total variation distance is to use Hellinger distance, defined
for $\nu = h \mu$ by
\begin{equation}\label{eqhellinger}
d_H(\nu,\mu) = 2 \, \int (1 - \sqrt h) \, d\mu \, .
\end{equation}
It is elementary to check that
\begin{equation}\label{eqcompare}
d_H(\nu,\mu) \, \leq \, 2 \, \parallel \mu - \nu \parallel_{TV} \, \leq 4 \, \sqrt{d_H(\nu,\mu)}
\end{equation}
hence both distances are ``almost'' equivalent. Using the concavity of $u \mapsto \sqrt{u}$ it is
also immediate that
\begin{equation}\label{eqhell2}
d_H\left(\frac{\nu+\mu}{2},\mu\right) \leq \frac 12 \, d_H(\nu,\mu) \leq \parallel \mu - \nu
\parallel_{TV} = 2\,\parallel \mu - \frac{\nu+\mu}{2} \parallel_{TV} \leq 8 \,
\sqrt{d_H\left(\frac{\nu+\mu}{2},\mu\right)}
\end{equation}
so that (with some changes in the constants) we may assume that $\nu=h \mu$ with $h\geq 1/2$.
\smallskip

Introduce as usual $I(t)=d_H(P_t^*h \mu, \mu)$ for some density of probability $h$, and
differentiating w.r.t. $t$, we get
\begin{equation}\label{eqhell3}
\frac {d}{dt} \, I(t) = - \frac 14 \, \int \frac{|\nabla P_t^* h|^2}{(P_t^*h)^{3/2}} \, d\mu \, .
\end{equation}
As in the preceding subsections we may state

\begin{proposition}\label{prophellinger}
Assume that there exists some non-increasing function $\beta_{H}$ defined on $(0,+\infty)$ such
that for all $s>0$ and all $f$ belonging to $D_2(L)$ the following inequality holds
\begin{equation}\label{eqineqhellinger}
\left(\int \, f^4 \, d\mu \right)^{1/2} \, - \, \int f^2 d\mu \leq \, \beta_{H}(s) \, \int
\Gamma(f) d\mu + s \, \Osc(f^2) \, ,
\end{equation}
then for all $\nu = h \mu$ , $d_H(P_t^*h \mu,\mu) \leq 3 \xi_H(t) \, \parallel h
\parallel_\infty^{1/2}$ with
$$\xi_{H}(t) = \inf \, \{s>0 \, , \, \beta_{H}(s) \, \log (1/s) \, \leq \, 4t \} \, .$$ Hence, if
$\tilde{\eta}(u) = u^{1/4} \phi(u)$,
$$\parallel P_t^*\nu - \mu \parallel_{TV} \leq \, \frac{4 \, \int h \phi(h) d\mu}{(\phi \circ \tilde{\eta}^{-1})
\left( 2 \, \int h \phi(h) d\mu/ \sqrt{3 \xi_H(t)}\right)} \, .$$
\end{proposition}
\begin{proof}
Apply \eqref{eqineqhellinger} with $f=(P_t^* h)^{1/4}$. It yields $$\frac {d}{dt} \, I(t) \leq -
\frac{4}{\beta_H(s)} I(t) + \frac{4s}{\beta_H(s)} \parallel P_t^* h \parallel_\infty$$ hence the
result (because $\Osc(h^{1/2}) \leq \parallel h \parallel_\infty^{1/2}$ and $d_H(\nu,\mu) \leq
2$).
\end{proof}

Note that \eqref{eqineqhellinger} implies the following
\begin{equation}\label{eqineqhellingervar}
\Var_\mu(f^2) \, \leq \, 2 \, \left(\beta_{H}(s) \, \int \Gamma(f) d\mu + s \, \Osc(f^2)\right) \,
\left(\int f^4 d\mu\right)^{1/2} \, ,
\end{equation}
just multiplying both hand sides in \eqref{eqineqhellinger} by $\left(\int f^4 d\mu\right)^{1/2} +
\int f^2 d\mu$ and applying Cauchy-Schwarz inequality. Conversely \eqref{eqineqhellingervar}
implies \eqref{eqineqhellinger} up to a factor 2 (majorizing $(\int f^4 d\mu)^{1/2}$ in the right
hand side by $\left(\int f^4 d\mu\right)^{1/2} + \int f^2 d\mu$ and then dividing both hand sides
by this quantity).
\medskip

Using $f=1+\varepsilon g$ for $\varepsilon$ going to 0 and $g$ bounded, we immediately see that
\eqref{eqineqhellingervar} implies a weak Poincar\'e inequality with $\beta_{WP}(s)=\frac 12 \,
\beta_H(2s)$. But this result can be greatly improved as follows.

\begin{proposition}\label{prophellingercapacity1}
If $\mu$ satisfies \eqref{eqineqhellinger} then for all $A$ s.t. $0 < \mu(A) < 1/2$ , $$Cap_\mu(A)
\, \geq \, \left(\frac{\sqrt 2 - 1}{2 \sqrt 2}\right) \,
\left(\frac{\mu^{1/2}(A)}{\beta_H\left(\frac{\sqrt 2 - 1}{2 \sqrt 2} \,
\mu^{1/2}(A)\right)}\right) \, . $$ Conversely if $Cap_\mu(A) \, \geq \,
\frac{\mu(A)}{\gamma(\mu(A))}$ for some non-increasing positive function $\gamma$, then
\eqref{eqineqhellingervar} holds with $2 \, \beta_H(s) = 48 \, \frac{\gamma(s^2)}{s}$ .
\end{proposition}
\begin{proof}
We start with the proof of the direct part. Let $\BBone_A \leq f \leq \BBone_\Omega$ with
$\mu(\Omega) \leq 1/2$. Then
$$\int f^2 d\mu = \int \BBone_\Omega \, f^2 d\mu \leq \left(\mu(\Omega)\right)^{1/2} \, \left(\int
f^4 d\mu\right)^{1/2} \leq (1/\sqrt 2) \, \left(\int f^4 d\mu\right)^{1/2} \, ,$$ so that
$$\left(\int f^4 d\mu\right)^{1/2} -  \int f^2 d\mu \geq (\sqrt 2 - 1/\sqrt 2) \, \left(\int f^4
d\mu\right)^{1/2} \geq (\sqrt 2 - 1/\sqrt 2) \, \mu^{1/2}(A) \, .$$ Since $0 \leq f \leq 1$,
$\Osc(f^2) \leq 1$. The result follows from \eqref{eqineqhellinger} with $s=\frac{\sqrt 2 - 1}{2
\sqrt 2} \, \mu^{1/2}(A)$.

For the converse part we use \eqref{eqineqhellingervar} and the proof of Theorem 2.2 in
\cite{BCR2} as modified in \cite{DGGW} Theorem 5.3 and Lemma 5.2. Indeed both Theorems are written
for the usual $\Gamma(f)=|\nabla f|^2$ on a riemanian manifold but $\mu$ absolutely continuous
with respect to the volume measure in \cite{BCR2} while this assumption is skipped in \cite{DGGW}.
The latter can be extended to our framework without any change.

By homogeneity we may assume that $\int f^4 d\mu =1$. In order to control
$\Var_\mu(f^2)$ we introduce a median $m$ of $f^2$ and use as usual $$\Var_\mu(f^2) \leq
\int_{\Omega_+} \, (f^2 - m)_+^2 d\mu + \int_{\Omega_-} \, (f^2 - m)_-^2 d\mu$$ with
$\Omega_+=\{f^2 \geq m\}$ and $\Omega_-=\{f^2 \leq m\}$. Define $g=(f^2-m)_+$. For a given $s>0$
we may choose $c=c(s):=\inf\{u\geq 0 \, ; \, \mu(g>u)\leq s\}$ so that $\mu(g>c) \leq s$. The case
$c=0$ is similar to \cite{BCR2,DGGW}. So we assume that $c>0$. It holds
$$\int_{g>c} \, g^2 \, d\mu \, \leq \, \parallel g \parallel_\infty \, \int_{g>c} \, g \, d\mu \,
\leq \, \parallel g \parallel_\infty \, \sqrt{\mu(g>c)} \, \left(\int \, g^2 \,
d\mu\right)^{1/2}$$ but the latter is less than $(\int f^4 d\mu)^{1/2} =1$. So $$\int_{g>c} \, g^2
\, d\mu \, \leq \, \sqrt s \, \parallel g \parallel_\infty \, \leq \, \sqrt s \, \Osc(f^2) \, .$$
Now we may follow the proof of \cite{BCR2,DGGW} and introduce the level sets $\Omega_k=\{g> c
\rho^k\}$ for $0<\rho< 1$ and $k \in \N$. The only difference is that we have to compare
$\int_{\Omega_{k+1}\backslash \Omega_k} \Gamma(g) \, d\mu$ with $\int_{\Omega_{k+1}\backslash
\Omega_k} \Gamma(f) \, d\mu$ in order to obtain \eqref{eqineqhellingervar}. Since $\Gamma(g) \leq
4 \, f^2 \, \Gamma(f)$, we have $$\int_{\Omega_{k+1}\backslash \Omega_k} \Gamma(g) \, d\mu \, \leq
\, 4 \, (m+c \rho^k) \, \int_{\Omega_{k+1}\backslash \Omega_k} \Gamma(f) \, d\mu \, .$$ But thanks
to Markov inequality and since $\int f^4 d\mu =1$ for $\varepsilon >0$ , $$ s \, \leq \, \mu(g>
(1-\varepsilon) c) \, = \mu(f^2> m + (1-\varepsilon) c) \, \leq \, \frac{1}{\left(m +
(1-\varepsilon) c\right)^2}$$ so that $m + (1-\varepsilon) c \leq \sqrt{1/s}$. For $k\geq 1$ we
thus have $m+c \rho^k \leq \sqrt{1/s}$. For $k=0$ since $\int f^4 d\mu =1$ again we know that
$m\leq \sqrt 2$, so that for $s$ small enough the previous inequality is satisfied. Arguing as in
\cite{BCR2,DGGW} we thus have obtained $$\int \, g^2 \, d\mu \, \leq \, \sqrt s \, \Osc(f^2) \, +
\, \frac{4 (1+\rho) \, \gamma(s)}{\rho^2 \, (1-\rho) \, \sqrt s} \, \int_{\Omega_+} \, \Gamma(f)
\, d\mu \, .$$ The case of $\Omega_-$ is similar and easier. Indeed on $\Omega_-$, $f^2$ is
bounded by $m$ hence by $\sqrt 2$ so that we obtain a better inequality. But since we have to sum
up both, this is not relevant. The result follows for $\rho=1/2$.
\end{proof}

\begin{corollary}\label{corhellcap1}
Define $\gamma_H(s)= s^{1/2} \, \beta_H\left(\frac{\sqrt 2 - 1}{2 \sqrt 2} \, s^{1/2}\right)$ .
\begin{itemize}
\item  If $\mu$ satisfies \eqref{eqineqhellinger} and $s \mapsto \gamma_H(s)$ is non-increasing on
$(0,1/2)$, $\mu$ satisfies a weak Poincar\'e inequality with $\beta_{WP}(s)= 12 \, \gamma_H(s)$,
and conversely this weak Poincar\'e inequality implies \eqref{eqineqhellinger} for $\beta(s)=c
\beta_H(c's)$ where $c$ and $c'$ are some universal constants.\item If $\mu$ satisfies
\eqref{eqineqhellinger}, $s \mapsto \gamma_H(1/s)=\theta_H(s)$ is non-increasing and $s \mapsto s
\theta_H(s)$ is non-decreasing on $(2,+\infty)$, $\mu$ satisfies a super-Poincar\'e inequality
\eqref{eqsuperpoinc} with $\beta_{SP}(s)=8 \gamma_H(1/s)$ for $s\geq 2$ and $\beta_{SP}(s)=8
\gamma_H(1/2)$ for $1\leq s\leq 2$.
\end{itemize}
\end{corollary}
The result follows from the previous Theorem, \cite{BCR2} Theorem 2.2 and Theorem 2.1,  and
\cite{BCR3} Corollary 6. Actually both Theorems are written for the usual $\Gamma(f)=|\nabla f|^2$
on a riemanian manifold and $\mu$ absolutely continuous with respect to the volume measure. A
careful reading shows that Theorem 1, and hence Corollary 6 in \cite{BCR3} can be extended to our
general framework (the final argument in the proof of the aforementioned Theorem is not
necessary). We already discussed the case of \cite{BCR2} Theorem 2.2.
\smallskip

The second part of Corollary \ref{corhellcap1} can be improved thanks to the results in
\cite{BCR1}. Indeed since \eqref{eqineqhellinger} is equivalent (up to some constants) to a
capacity-measure criterion, it is equivalent to a general Beckner-type inequality (see \cite{BCR1}
section 5.3 for the definitions and Theorem 18 for the result).

In particular if $\beta_H(s)=c/s$, \eqref{eqineqhellinger} is equivalent to a Poincar\'e
inequality, and if $\beta_H(s)=c/s \log(1/s)$ it is equivalent to a logarithmic Sobolev
inequality. Notice that in the first case a direct application of Proposition \ref{prophellinger}
for $h$ such that $\int h^2 d\mu < +\infty$, i.e. with $\phi(u)=u$, yields a polynomial decay $c/
t^{2/5}$ which is disastrous, since Poincar\'e inequality yields an exponential decay.
\smallskip

Now, \eqref{eqineqhellinger} with $\beta_H$ constant is equivalent to the exponential decay $I(t)
\leq e^{- \alpha t} \, I(0)$ for some $\alpha > 0$, which implies according to \eqref{eqcompare},
$\parallel P_t^*h - 1\parallel_{\L^1(\mu)} \leq 2 \sqrt 2 \, e^{- \alpha t/2}$. But this implies a
super Poincar\'e inequality with $\beta_{SP}= c \, s^{-1/2}$, hence again $P_t$ is
ultracontractive according to Wang's result.

Hence, the direct study of the Hellinger distance furnishes no convincing results. However,  in
\cite{DGGW}, where inequalities in the spirit of (\ref{eqineqhellingervar}) were introduced under
the name of $L^q$-Poincar\'e inequalities, applications of those type of inequalities concern
large time behavior of nonlinear diffusions, namely porous media equation $\partial_tu=L(u^m)$ for
$m\ge 1$. Formal calculations indicate that (\ref{eqineqhellingervar}) could have the same role
for other nonlinear diffusions. We leave this for further research.
\medskip

%%%%%%%%%%%%%%%%%%%%%%%%%%%%%%%
%%%%%%%%%%%%%%%%%%%%%%%%%%%%%%%
%%%%%%%%%%%%%%%%%%%%%%%%%%%%%%%

\section{\bf Other related inequalities, reversing the roles.}\label{secother}

One of the main feature of the use of functional inequalities for studying the total variation
distance, is that symmetry is broken. Indeed $I_\psi(\Q|\P)$ is in general not symmetric. If it
seems natural to privilege the invariant measure $\mu$ by looking at $I_\psi(P_t^* \nu|\mu)$, one
may ask what happens if we reverse the roles. This idea is not completely new since in \cite{DLM}
the authors have studied the evolution of the total variation distance between $P_t^* \nu$ and
$P_t^* \nu'$ for any initial $\nu$ and $\nu'$, but under strong conditions on one of them.

In second place, since
\begin{equation}\label{eqminore}
\parallel P_t^* h - 1\parallel_{\L^1(\mu)} = 2 \, \parallel P_t^* \left(\frac{h+1}{2}\right)
 - 1\parallel_{\L^1(\mu)}
\end{equation}
we may assume that $h \geq \frac 12$, i.e. $P_t^* h \geq \frac 12$. Thus if $\nu=h \mu$
$$\frac{d\mu}{d P_t^* \nu} = \frac{1}{P_t^* h} \leq 2 \, .$$

We thus have, denoting $P_t^*\nu=\nu_t$,
\begin{eqnarray*}
\parallel P_t^*\nu - \mu \parallel_{TV} & \leq & \sqrt{\Var_{\nu_t}(1/P_t^* h)} \\
\parallel P_t^*\nu - \mu \parallel_{TV} & \leq & \sqrt{2 \, \Ent_{\nu_t}(1/P_t^* h)}
\end{eqnarray*}
so that we shall study
\begin{equation}\label{eqvarreverse}
V(t) = \Var_{\nu_t}(1/P_t^* h) = \int \, \frac{1}{P_t^* h} d\mu - 1 \, ,
\end{equation}
and
\begin{equation}\label{eqentreverse}
E(t) = \Ent_{\nu_t}(1/P_t^* h) = \int \, \log(1/P_t^*h) \, d\mu \, .
\end{equation}
Assuming first that $h$ is also bounded from above (for the forthcoming calculation to be
rigorous), we immediately get using the chain rule
\begin{equation}\label{eqreversederive}
\frac{d}{dt} \ V(t) = - \int \, \frac{1}{(P_t^* h)^3} \, \Gamma(P_t^*h) \, d\mu \quad \textrm{ and
} \quad \frac{d}{dt} \ E(t) = - \frac 12 \, \int \, \frac{1}{(P_t^* h)^2} \, \Gamma(P_t^*h) \,
d\mu \, .
\end{equation}
Remark now that the exponential decay $$V(t) \, \leq
\, e^{-\lambda t} \, V(0)$$ is equivalent to
\begin{equation}\label{eqvt}
\int (1/P_t^* h) d\mu - 1 \, \leq \, (1/\lambda) \, \int \, \frac{1}{(P_t^* h)^3} \,
\Gamma(P_t^*h) \, d\mu
\end{equation}
for all $t \geq 0$, and that the exponential decay $$E(t) \, \leq \,
e^{-\lambda t} \, E(0)$$ is equivalent to
\begin{equation}\label{eqvtlog}
\int \log (1/P_t^* h) d\mu  \, \leq \, (1/2 \lambda) \, \int \, \frac{1}{(P_t^* h)^2} \,
\Gamma(P_t^*h) \, d\mu
\end{equation}
for all $t \geq 0$.
\smallskip

There are now two approaches which can be seen as static or dynamic: the static one is to consider
equations (\ref{eqvt}) and (\ref{eqvtlog}) as functional inequalities and as before look at
capacity-measure conditions for these inequalities ; the dynamic one starts from the assumption
that $hd\mu$ satisfies some inequalities (say Poincar\'e for example) and  study the propagation
along the semigroup of such an inequality which enables us to get a direct control of $V(t)$. As
we will see, the static approach seems to be very restrictive, whereas under curvature assumptions
the dynamic one furnishes interesting result.

\subsection{\bf The static approach.}\label{secreversevariance}
\subsubsection{\bf Variance control.}
 Using \eqref{eqvt} with $u=\sqrt{1/P_t^* h}$, such an exponential decay for
all $h$ ($\geq 1/2$) is equivalent to
\begin{equation}\label{eqwefort}
\int \, u^2 \, d\mu \, - \, 1  \, \leq \, C_{WE} \, \int \Gamma(u) d\mu
\end{equation}
for all  $u$ belonging to $D_2(L)$  such that $0 \leq u \leq \sqrt{2}$ and $\int (1/u^2) d\mu =
1$. The weak version
\begin{equation}\label{eqwe}
\int \, u^2 \, d\mu \, - \, 1  \, \leq \, \beta_{WE}(s) \, \int \Gamma(u) d\mu + s  \, ,
\end{equation}
for some non-increasing function $\beta_{WE}$ defined on $(0,+\infty)$ and all $s>0$ implies that
for all $\nu = h \mu$ ,
\begin{equation}\label{eqtralala}
\parallel P_t^*\nu - \mu
\parallel_{TV} \leq \sqrt{V(t)} \leq C \sqrt{\xi_{WE}(t)} \, ,
\end{equation}
for some universal constant $C$ where $\xi_{WE}(t) = \inf \, \{s>0 \, , \, \beta_{WE}(s) \, \log
(1/s) \, \leq \, 4t \}$.
\smallskip

If we relax the condition $u \leq \sqrt 2$, the inequalities (\ref{eqwefort}) and (\ref{eqwe}) are
extremely strong if $\mu$ has no atoms. Indeed, if $\int (1/u^2) d\mu < +\infty$, \eqref{eqwe}
becomes
$$ \int \, u^2 \, d\mu \, \leq \, \beta_{WE}(s) \, \int \Gamma(u) d\mu + (1+s)  \, \frac{1}{\int
(1/u^2) d\mu} \, .$$ Let $u$ be such that $essinf (u) =0$, so that, since $\mu$ has no atoms, for
all $\varepsilon >0$, $\mu(u \leq \varepsilon) > 0$. Choose $f=(u - \varepsilon)_+ + \chi^2$ and
apply \eqref{eqwe}. It yields $$ \int \, f^2 \, d\mu \, \leq \, \beta_{WE}(s) \, \int_{u\geq
\varepsilon} \Gamma(u) d\mu + (1+s)  \, \frac{1}{\int (1/f^2) d\mu} \, .$$ Now we may let $\chi$
go to 0, so that we obtain $$ \int \, (u - \varepsilon)_+^2 \, d\mu \, \leq \, \beta_{WE}(s) \,
\int \Gamma(u) d\mu$$ for all $s$, so that we may replace $\beta_{WE}(s)$ by
$\beta_{WE}(1)=\beta$. Next we let $\varepsilon$ go to 0 and obtain for all $u$ such that $essinf
(u) =0$, $$\int \, u^2 \, d\mu \, \leq \, \beta \, \int \, \Gamma(u) \, d\mu \, .$$ If $\BBone_A
\leq f \leq \BBone_\Omega$ as usual, applying the previous inequality with $u=1-f$ and $\mu(A)>0$
so that $essinf (u) =0$) yields $1/2 \leq \beta \, Cap_\mu(A)$ for all $A$, in particular of
course ultracontractivity.
\smallskip

This discussion indicates that if we stay with $u \leq \sqrt 2$, a natural choice is $u = \sqrt 2
\, (1 - \alpha f)$ for some $0 \leq \alpha \leq 1$. Note that for $\alpha=0$, $\int (1/u^2) d\mu =
1/2$ while for $\alpha =1$ it is equal to $+\infty$ as soon as $\mu(A)>0$. By monotonicity and
continuity we may thus find a unique $\alpha_0$ such that $\int (1/u^2) d\mu = 1$. Hence
\begin{eqnarray*}
1  & = & \int_{\Omega^c} (1/u^2) d\mu + \mu(\Omega) \int_\Omega (1/u^2) \,
\frac{d\mu}{\mu(\Omega)} \\ & = & \frac 12 \, \mu(\Omega^c) + \mu(\Omega) \int_\Omega (1/u^2) \,
\frac{d\mu}{\mu(\Omega)} \\ & \geq & \frac 14 + \mu(\Omega) \int_\Omega (1/u^2) \,
\frac{d\mu}{\mu(\Omega)} \, ,
\end{eqnarray*}
so that using Jensen inequality $\int (1/u^2) d\nu \geq 1/\left( \int u^2 d\nu\right)$ we obtain
$$\int_\Omega \, u^2 \, \frac{d\mu}{\mu(\Omega)} \geq \frac{4 \, \mu(\Omega)}{3} \, .$$ \eqref{eqwe}
thus implies $$2 \mu(\Omega^c) + \frac{4 \, \mu^2(\Omega)}{3} - 1 = 1 - 2 \mu(\Omega) + \frac{4 \,
\mu^2(\Omega)}{3} \leq 2 \, \alpha_0^2 \, \beta_{WE}(s) \int \Gamma(f) d\mu \, + \, s \, .$$ But
the minimum of $1 - 2x + (4/3) x^2$ is attained for $x=3/4$ and is equal to $1/4$. Choosing
$s=1/8$ for instance and since $\alpha_0 < 1$, we thus have shown that there exists $\theta>0$
such that $Cap_\mu(A) \geq \theta$ for all $A$ with $\mu(A)>0$. We are thus in the
ultracontractive situation, i.e. we have a uniform exponential decay in $\L^1$ and not the poor
one given by \eqref{eqtralala}.

\begin{remark}\label{remcappoint}
Note that if $A$ is non empty but $\mu(A)=0$, we may find some $\Omega$ with $\mu(\Omega) \leq
1/2$, some $f$ such that $\BBone_A \leq f \leq \BBone_\Omega$ so that $Cap_\mu(A) \geq \, \int
\Gamma(f) d\mu - (\theta/2)$. We may assume that $f$ is uniformly continuous according to our
assumptions on the model. If for all $1>c>0$ , $\mu(f>c)>0$ (for instance if $\mu(B)>0$ for any
non empty open ball $B$) define $g=\min (1 ; f/c)$. Then $\BBone_C \leq g \leq \BBone_\Omega$ for
$C=\{f\geq c\}$ so that $$\theta \leq Cap_\mu(C) \leq \int \Gamma(g) d\mu \leq (1/c^2) \int
\Gamma(f) d\mu \leq (1/c^2) \, \left(Cap_\mu(A) + \frac{\theta}{2}\right)$$ so that $Cap_\mu(A)
\geq (c^2 - (1/2)) \, \theta$ for all $c$ hence $Cap_\mu(A) \geq \theta/2$. \hfill $\diamondsuit$
\end{remark}
\medskip

\subsubsection{\bf Entropy control.}\label{secreverseentropy}

We focus now on the analysis of \eqref{eqvtlog}. Here again we may state :

assume that there exists some non-increasing function $\beta_{MT}$ defined on $(0,+\infty)$ such
that for all $s>0$ and all $v$ belonging to $D_2(L)$ such that $v \geq - \log 2$ and $\int v d\mu
=0$, the following inequality holds
\begin{equation}\label{eqwmt}
\log \, \left( \int e^v \, d\mu\right) \, \leq \, \beta_{MT}(s) \, \int \Gamma(v) d\mu + s \,
\Osc^2(v) \, ,
\end{equation}
then for all $\nu = h \mu$ , $$\parallel P_t^*\nu - \mu \parallel_{TV} \leq  \sqrt{2 \,
\xi_{MT}(t)} \, \sqrt{\log 2 + \Osc^2(\log h)} \, ,$$ for some universal constant $C$, where
$$\xi_{MT}(t) = \inf \, \{s>0 \, , \, 2 \beta_{MT}(s) \, \log (1/s) \, \leq \, t \} \, .$$
Hence for all $\nu = h \mu$,  for $\eta(u)= \phi(u) \, \log(u)$, $$\parallel P_t^*\nu - \mu
\parallel_{TV} \leq \frac{4 \, \int h \, \phi(h) \, d\mu}{(\phi \circ \eta^{-1})\left( \left(\int
h \, \phi(h) \, d\mu\right)/\sqrt{\xi_{MT}(t)}\right)} \, .$$ Indeed, if \eqref{eqwmt} holds, we
may choose
$$v_t = \log(P_t^* h) - \int \log(P_t^* h) d\mu \, ,$$ and apply \eqref{eqreversederive}. We
obtain $$\frac{d}{dt} \ E(t) \leq \frac{-1}{2\beta_{MT}(s)} \, E(t) \, + \, \frac{s}{2
\beta_{MT}(s)} \, \Osc^2(v_t) \, .$$ Gronwall's lemma immediately yields $E(t) \leq \xi_{MT}(t)
(E(0)+\Osc^2(\log h))$ because $\Osc^2(\log P^*_t h) \leq \Osc^2(\log h)$. Since $E(0) \leq \log
2$ if $h\geq \frac 12$ the conclusion follows.

Inequality \eqref{eqwmt} is a weak version of the so called Moser-Trudinger inequality, i.e. for
all nice $v$ such that $\int v d\mu =0$, $$\log \, \left( \int e^v \, d\mu\right) \, \leq \,
C_{MT} \, \int \Gamma(v) d\mu \, ,$$ for some constant $C_{MT}$, which appears as some limit case
of Sobolev inequalities. As for equation \eqref{eqwe} we shall see that it implies a very strong
capacity measure inequality.
\smallskip

First \eqref{eqwmt} is equivalent to $$\int \log(1/u) d\mu \, \leq \, \beta_{MT}(s) \, \int
\frac{\Gamma(u)}{u^2} \, d\mu \, + \, s \, \Osc^2(\log(u)) \, ,$$ for $1/2 \leq u$ and $\int u
d\mu =1$. With $v=(1/u^2)$ so that $0\leq v \leq 4$ and $\int (1/\sqrt v) d\mu =1$ we obtain
$$\int \, \log(v) \, d\mu \, \leq \, \frac 12 \, \beta_{MT}(s) \, \int
\Gamma(v) \, d\mu \, + \, \frac s4 \, \Osc^2(\log(v)) \, .$$ For $\BBone_A \leq f \leq
\BBone_\Omega$ with $\mu(\Omega)\leq 1/2$ choose $v =  4 \, (1 - \alpha f)$ for some $0 \leq
\alpha \leq 1$. Then $0 \leq v \leq 4$ and $\int (1/\sqrt v) d\mu$ is equal to $1/2$ for
$\alpha=0$ and goes to $+\infty$ when $\alpha$ goes to $1$, provided $\mu(A)>0$. We thus may
choose $\alpha_1$ such that this integral is equal to $1$.

Now $\int_\Omega (1/\sqrt v) d\mu = 1 - (1/2) \mu(\Omega^c) \leq 3/4$. Since $x \log(x) \geq -
1/e$ for $0<x \leq 1$, it implies
\begin{eqnarray*}
\int \, \log(v) \, d\mu & \geq & \log(4) \, \mu(\Omega^c) + 2 \int_{v\leq 1} \, \log(\sqrt v) \,
d\mu \\ & \geq & \log(4) \, \mu(\Omega^c) - \frac 2e \, \int_{v \leq 1} \, (1/\sqrt v) d\mu \\
 & \geq & \log(4) \, \mu(\Omega^c) - \frac 2e \, \int_\Omega \, (1/\sqrt v) d\mu \\
& \geq & \log(4) \, \mu(\Omega^c) - \frac 2e \, (1 - (1/2) \mu(\Omega^c)) \\ & \geq & (\log(4) +
(1/e)) \, \mu(\Omega^c) - (2/e) \, \geq \, \frac 12 (\log(4) - (3/e)) \, \geq \, 0.1 \, .
\end{eqnarray*}
Again this implies (choosing $s=0.05$) that $Cap_\mu(A) \, \geq \, c > 0$ for all $A$ such that
$\mu(A) > 0$, hence ultracontractivity.

\begin{remark}\label{remMT}
We already saw in Remark \ref{remcappoint} that under mild assumptions ($\mu(B)>0$ for all non
empty open ball for instance), we may deduce that $Cap_\mu(A) \geq c/2$ for all non empty set $A$.

If $\mu$ is an absolutely continuous measure with respect to the riemanian measure on a riemanian
manifold $M$, with an everywhere positive density, this assumption is satisfied. If $M$ is non
compact, it is easily seen that one can build Lipschitz function $h$ vanishing on balls $B(x_0,R)$
such that $\parallel h \parallel_{Lip} \leq 1/R$ so that the capacity of points with $d(x,x_0)>
3R$ has to be smaller that $1/R^2$. Hence for a non compact riemanian manifold, the capacity of
all points cannot be bounded from below. It follows that the Moser-Trudinger inequality (even in
its weak form) cannot be satisfied on non compact manifolds. \hfill $\diamondsuit$
\end{remark}

\subsection{\bf The dynamic approach.}

As previously said, another direct approach of \eqref{eqvt}, close to \cite{DLM}, is the following.

Assume that $\nu_t=P_t^*h \, \mu$ satisfies a Poincar\'e inequality
\begin{equation}\label{eqpoincnut}
\Var_{\nu_t}(g) \, \leq \, C_P(t) \, \int \, \Gamma(g) \, d\nu_t \, .
\end{equation}
Applying \eqref{eqpoincnut} with $g=1/P_t^*h$ yields precisely \eqref{eqvt} with $\lambda =
1/C_P(t)$ and consequently
\begin{equation}\label{eqvtah}
V(t) \, \leq \, e^{ - \, \int_0^t \, (1/C_P(s)) \, ds} \, V(0) \, .
\end{equation}
Here we have to be much more accurate. Indeed in the derivation of \eqref{eqtralala} we may first
establish \eqref{eqwe} and \eqref{eqreversederive} for nice $h$'s bounded from below and from
above, and then extend \eqref{eqtralala} to general $h$'s. Here, since we are using Poincar\'e
inequality for $\nu_t$, we need \eqref{eqreversederive} for the $h$ we are interested in. It is
not difficult to see that $h \in \L^2(\mu)$ and $h\geq 1/2$ are sufficient for all these
derivations to be correct. Dealing with initial densities in $\L^2(\mu)$ is certainly
disappointing, however we shall see how to use approximations by such functions, but to this end
we have to weaken our assumptions.
\smallskip

From now on we assume that $1/2 \leq h \leq K$ is a (nice) bounded density of probability such
that $\nu=h \mu$ satisfies a weak Poincar\'e inequality with function $\beta_{WP}$, we expect not
depending on $K$. \eqref{eqreversederive} is thus satisfied, and we want to study a possible weak
Poincar\'e inequality for $\nu_t=P_t^*h \mu$. It holds
\begin{eqnarray*}
\Var_{\nu_t}(f) & = & \int f^2 \, P_t^*h \, d\mu \, - \, \left(\int \, f \, P_t^*h \,
d\mu\right)^2 \, = \, \int P_t(f^2) \, h \, d\mu \, - \, \left(\int \, P_t f \, h \,
d\mu\right)^2 \\ & = & \int P_t(f^2) \, h \, d\mu \, - \, \int \left(P_t f\right)^2 \, h \, d\mu
\, + \, \int \left(P_t f\right)^2 \, h \, d\mu \, - \, \left(\int \, P_t f \, h \, d\mu\right)^2\\
& \leq & \left(\int     P_t(f^2) \, h \, d\mu \, - \, \int \left(P_t f\right)^2 \, h \, d\mu\right)
\, + \, \beta_{WP}(s) \, \int \Gamma(P_t f) \, h \, d\mu  \, + \, s \, \Osc^2(f) \, .
\end{eqnarray*}
It remains to exchange $\Gamma$ and $P_t$ and to control the first term in the left hand side.
Both controls are known to be equivalent to a curvature assumption introduced by Bakry and Emery.
Recall some definitions (see e.g. \cite{logsob} section 5.3 and Proposition 5.4.1)
\begin{proposition}\label{defgamma2}
Introduce $\Gamma_2(f,g) := (1/2) \left(L \Gamma(f,g) - \Gamma(Lf,g) - \Gamma(f,Lg)\right)$. We
shall say that the curvature of $L$ is bounded from below by $\rho \in \R$ if for all nice $f$,
$\Gamma_2(f) := \Gamma_2(f,f) \geq \rho \, \Gamma(f)$.

Then the following three assertions are equivalent
\begin{itemize}
\item the curvature of $L$ is bounded from below by $\rho \in \R$, \item for all nice $f$ and all
$t>0$, $\Gamma(P_t f) \, \leq \, e^{-\rho t} \, P_t(\Gamma(f))$, \item for all nice $f$ and all
$t>0$, $P_t(f^2) - (P_t f)^2 \, \leq \, \frac{1 - e^{-\rho t}}{\rho} \, P_t(\Gamma(f))$ (for
$\rho=0$ replace the coefficient by $t$).
\end{itemize}
\end{proposition}
We immediately deduce

\begin{corollary}\label{corgamma2}
If the curvature of $L$ is bounded from below by $\rho \in \R$, and $\nu=h \mu$ satisfies a weak
Poincar\'e inequality with function $\beta_{WP}$, then $\nu_t = P_t^*h \mu$ satisfies a weak
Poincar\'e inequality with function $$\beta_{WP}(t,s)= \frac{1 - e^{-\rho t}}{\rho} \, + \,
e^{-\rho t} \, \beta_{WP}(s) \, .$$
\end{corollary}

Plugging this estimate (with $f= 1/P_t^*h$) in \eqref{eqreversederive} yields (since $\Osc(f) \leq
2$) $$\frac{d}{dt} \ V(t) \leq - \, \frac{V(t)}{\beta_{WP}(t,s)} + \frac{4s}{\beta_{WP}(t,s)} \,
.$$ Defining
\begin{equation}\label{eqratecurve}
r(t,s) = \int_0^t \, \frac{1}{\beta_{WP}(u,s)} \, du
\end{equation}
we obtain $$ V(t)  \leq  e^{-r(t,s)} \, V(0) + \int_0^t \, \frac{4s}{\beta_{WP}(u,s)} \,
e^{r(u,s)-r(t,s)} \, du  \,  \leq  \,  e^{-r(t,s)} \, V(0) \, + \, 4s \, .$$ But an elementary
calculation yields $$r(t,s)= \log \left(\frac{e^{\rho t} + \rho \, \beta_{WP}(s) -1}{\rho \,
\beta_{WP}(s)}\right) \, $$ so that finally $$V(t) \, \leq \, \frac{\rho \, \beta_{WP}(s)}{e^{\rho
t} + \rho \, \beta_{WP}(s) -1} + 4s \, .$$ Of course this inequality has some interest only if
$r(t,s) \to +\infty$ as $t \to +\infty$, hence if $\rho \geq 0$. It also yields
\begin{equation}\label{eqratecurve2}
\parallel P_t^*h \mu - \mu\parallel_{TV} \leq \left(\inf_{s>0} \left\{\frac{\rho \, \beta_{WP}(s)}{e^{\rho
t} + \rho \, \beta_{WP}(s) -1} + 4s \right\}\right)^{1/2} \, .
\end{equation}
This inequality is more tractable since it extends to any $h$ (up to a factor 2, recall
\eqref{eqminore}), if we can approximate $(h+1)/2$ by a sequence of $h_n$ with $n \geq h_n \geq
1/2$ such that each $h_n \mu$ satisfies a weak Poincar\'e inequality with the same $\beta_{WP}$.
As a consequence $((h+1)/2) \mu$ will satisfy the same inequality. We have thus shown

\begin{theorem}\label{thmdelmo}
Assume that the curvature of $L$ is bounded from below by $\rho \geq 0$. Assume that one can find
a sequence $h_n$ with $n \geq h_n \geq 1/2$ such that each $h_n \mu$ satisfies a weak Poincar\'e
inequality with the same $\beta_{WP}$, such that $h_n \to (1+h)/2$ in $\L^1(\mu)$ as $n$ goes to
infinity. Then \eqref{eqratecurve2} holds.

In particular if $s \mapsto \beta_{WP}(s)/s$ is non-increasing (at least for small $s$), define
$\theta(u) = \inf \, \{s \, ; \, (\beta_{WP}(s)/s) \leq 4u/\rho \}$. Then there exists a constant
$C$ such that $$\parallel P_t^*h \mu - \mu\parallel_{TV} \leq C \, \theta^{1/2}(e^{\rho t}) \, .$$
In particular if $\beta_{WP}(s) \leq c s^{-q}$ for some $q\geq 0$, $\parallel P_t^*h \mu -
\mu\parallel_{TV} \leq C \, e^{- \rho t/2(1+q)}$.
\end{theorem}

Recall that $\rho>0$ implies that $\mu$ satisfies a log-Sobolev inequality with $C_{LS}=2/\rho$.
Thus, since $h_n$ is bounded below (by $1/2$) and above (say by $n$), $h_n \mu$ satisfies a
log-Sobolev inequality with a constant depending on $n$. Using this constant and estimating
$\parallel h \mu - h_n \mu\parallel_{TV}$ is similar (actually a little bit worse) to the
truncation method we used in subsection \ref{seclogsob}.

Nevertheless, since $\mu$ satisfies a log-Sobolev inequality, we have to compare the result
obtained in \eqref{eqratecurve2} and the ones in Corollary \ref{corlogsob}, for densities $h$
which are not in the space $\L \log \L$ (nor in any $\L \log_+^{\beta} \L$ for $\beta>0$ according
to \eqref{eq4} in Example \ref{examplepentropy}). Since there is no general criterion for the weak
Poincar\'e inequality, we shall make this comparison on examples only.

\begin{example}\label{exdelmo}
Let us assume that $\mu(dx)=e^{-V(x)} dx$ is a probability measure on $\R$, $Lf= f'' - V'f'$ so
that $\mu$ is a symmetric measure for $L$ and $\Gamma(f)=2 (f')^2$. If $V''(x) \geq \rho
> 0$, then the curvature of $L$ is bounded from below by $\rho$.

Choose $h=e^V \, g$, with $g\geq 0$ and $\int g dx =1$. For simplicity we assume that $g$ is
symmetric. Then it is known (see \cite{BCR2} Theorem 2.3) that $\nu = h \mu$ satisfies a weak
Poincar\'e inequality with a function $\beta_{WP}(s)= C \, \beta(s)$, if $\beta$ is a non
increasing function, for $12 B \geq C \geq (1/4) b$ where $$b=\sup_{x>0} \,
\frac{\nu([x,+\infty))}{\beta\left(\frac{\nu([x,+\infty))}{4}\right)} \, \int_0^x (1/g)(y) dy \,
\, \textrm{ and } \, \, B=\sup_{x>0} \, \frac{\nu([x,+\infty))}{\beta(\nu([x,+\infty)))} \,
\int_0^x (1/g)(y) dy \, .$$ We immediately see that if $g \gg e^{-V}$, $((1+h)/2) \mu$ will
satisfy a weak Poincar\'e inequality with the same $\beta$ and a modified constant $C$, and that
$((1+h\wedge n)/2) \mu$ will also satisfy a weak Poincar\'e inequality with the same $\beta$ and a
constant $D$ which can be chosen independent of $n$.

As in Remark \ref{remlemeilleur} we may evaluate $$\int_0^x \, (1/g)(y) dy \sim \frac{-1}{g'(x)}
\quad \textrm{ and } \quad \int_x^{+\infty} g(y) dy \sim \frac{- g^2(x)}{g'(x)} \, ,$$ provided
$g'<0$ near infinity and $\lim_{x \to +\infty} (g(x) \, g''(x)/(g')^2(x)) = 1$ (see e.g.
\cite{logsob} Proposition 6.4.1).
\smallskip

We shall give some explicit examples
\begin{itemize}
\item If $\nu([x,+\infty)) \sim x^{-p}$ for some $p>0$, i.e. $g$ behaves like $x^{-(1+p)}$ at
infinity, $\nu$ satisfies a weak Poincar\'e inequality with $\beta_{WP}(s) \sim s^{-2/p}$, so that
Theorem \ref{thmdelmo} furnishes an exponential decay $e^{-\rho p t/2(p+2)}$. For such a result
using Corollary \ref{corlogsob} we need that $h \in \L \log_+^\beta L$ for some $\beta>0$, that is
we need $\int V^{\beta}(x) g(x) dx < +\infty$. Hence if $V(x) \sim x^k$ near infinity, we need
$\beta< p/k$ and obtain a decay slightly worse than $e^{-\rho p t/2(p+k)}$. In all cases (since
$k\geq 2$) this is a worse decay than $e^{-\rho p t/2(p+2)}$.

If $V \sim e^x$ at infinity, the situation is still worse since Corollary \ref{corlogsob} does not
furnish the exponential decay which is still true according to Theorem \ref{thmdelmo}. \item If
$g(x) \sim (1/ x \, \log^2(x))$ at infinity, we get $\beta_{WP}(u) \sim e^{2/\sqrt u}$. This
yields a polynomial decay $c/t$ for the total variation distance. Now it is easily seen that, if
$V(x)=x^2$, $\int h \log_+^{1-\varepsilon}(\log_+(h)) \, d\mu < +\infty$ for $\varepsilon >0$ and
infinite for $\varepsilon=0$. Corollary \ref{corlogsob} furnishes a decay $c/t^{(1 -
\varepsilon)/2}$, hence still a worse rate. Again for larger $V$'s the result is unchanged with
Theorem \ref{thmdelmo} and is getting worse with Corollary \ref{corlogsob}.
\end{itemize}

It seems in the one dimensional case that Theorem \ref{thmdelmo} gives better results than
Corollary \ref{corlogsob}, in particular because it does not take into account the moments of $h$
with respect to $\mu$ (this is not completely true since these moments have an influence on
$\beta_{WP}$ but we may change $\mu$ without changing nor $\rho$ nor $\beta_{WP}$).

A better understanding of general weak Poincar\'e inequalities is however necessary to claim that
it has to be a general fact. In addition, $\mu$ is supposed to satisfy a strong form of the
log-Sobolev inequality (the Bakry-Emery condition). Finally it is quite possible that for very
oscillating densities (not satisfying the conditions for the tail estimates to be true for
instance) one can have finite entropy but a bad weak Poincar\'e function. \hfill $\diamondsuit$
\end{example}
\medskip

\begin{remark}
We have previously studied the propagation of weak Poincar\'e inequalities along the semi-group.
It is then natural to look at the propagation of Super-Poincar\'e inequalities. Indeed, assume
that $\nu=h\mu$ satisfies a Super-Poincar\'e inequality (see \cite{BCR3} or \cite{BLW} for
explicit conditions on $h$ and $\mu$), i.e. there exists $\beta_{SP}$ defined on $[1,\infty[$ such
that
$$\Var_{\nu}(f)\le\beta_{SP}(s)\int \Gamma(f)d\nu+s\left(\int|f|d\nu\right)^2$$ and if the curvature
of $L$ is bounded below by $\rho$, then, by the same proof as before $\nu_t=P^*_th\mu$ satisfies a
Super-Poincar\'e inequality with $\beta_{SP}(t,s)=\rho^{-1}(1-e^{-\rho t})+e^{-\rho
t}\beta_{SP}(s)$. We can then as before (with the same precautions on $h$ as before) plug these
estimate into (\ref{eqreversederive}) to get
\begin{equation}\label{eqratecurve2sp}
\parallel P_t^*h \mu - \mu\parallel_{TV} \leq \left(\inf_{s>1} \left\{     \frac{\rho \, \beta_{WP}(s)}{e^{\rho
t} + \rho \, \beta_{WP}(s) -1} + (s-1)  \right\}\right)^{1/2} \, ,
\end{equation}
which does not give better result than a Poincar\'e inequality.\hfill$\diamondsuit$
\end{remark}

\smallskip

\begin{remark}
Concerning the same approach for entropy, let us point out the following remarks. First recall
Proposition 5.4.5 in \cite{logsob} which applies here since we are in the framework called ``the
diffusion case'' therein. Hence under the assumptions previously set we may replace the weak
Poincar\'e inequality by a weak log-Sobolev inequality, and (up to some constants) replace
$\beta_{WP}$ by $\beta_{WLS}$ in \eqref{eqratecurve2}. This is of course not very clever since
$\beta_{WLS}$ is much bigger than $\beta_{WP}$.\hfill$\diamondsuit$

\end{remark}
\medskip

\medskip

\bigskip
\bibliographystyle{plain}

%********************************************************************************

\end{document}